\documentclass[11pt]{article}
\usepackage{mathrsfs}
\usepackage{amscd}
\usepackage{amsmath}
\usepackage{latexsym}
\usepackage{amsfonts}
\usepackage{amssymb}
\usepackage{amsthm}
\usepackage{graphicx}
\usepackage{color}
\usepackage{cite}

\newtheorem{proposition}{Proposition}[section]

\newtheorem{lemma}{Lemma}[section]
\newtheorem{theorem}{Theorem}[section]

\theoremstyle{remark}
\newtheorem{remark}[theorem]{Remark}

\begin{document}
\title{On the Cauchy problem for the
 magnetic Zakharov system}
\author
         {
           {Boling Guo$^1$\ \ \ Jingjun Zhang$^{2,3}$\thanks{Corresponding author: zjj\_math@yahoo.com.cn }\ \ \ Chunxiao Guo$^{2}$ }\\
           {\small $^1$Institute of Applied Physics and Computational Mathematics}\\
           {\small P. O. Box 8009,\ Beijing,\ China,\ 100088}\\
           {\small $^2$The Graduate School of China Academy of Engineering Physics  }\\
           {\small P. O. Box 2101,\ Beijing,\ China,\ 100088 } \\
           {\small $^3$College of Mathematics and Information Engineering}\\
           {\small   Jiaxing University, Zhejiang, \ China,\ 314001 } \\
         \date{}
         }
\maketitle
\begin{abstract}
In this paper, we study the Cauchy problem of the magnetic type
Zakharov system which describes the pondermotive force and
magnetic field generation effects resulting from the non-linear
interaction between plasma-wave and particles. By using the energy
method to derive a priori bounds and an approximation argument for
the construction of solutions, we obtain local existence and
uniqueness results for the magnetic Zakharov system in the case of
$d=2,\ 3$.
\end{abstract}
\begin{center}
 \begin{minipage}{120mm}
\end{minipage}
\end{center}
\begin{center}
 \begin{minipage}{120mm}
\end{minipage}
\end{center}

%
%
%
%
%

\section{Introduction and the main results}
\setcounter{section}{1}\setcounter{equation}{0} In this paper we
study the Cauchy problem for the magnetic Zakharov system
\begin{equation}\label{1.1}
\left\{\!\!
\begin{array}{lc}
iE_{t}+\nabla(\nabla\cdot E)-\alpha\nabla\times(\nabla\times
E)-nE+i E\times B=0, &\\
n_{tt}-\triangle n=\Delta |E|^{2},&\\
B_{tt}+\Delta ^{2}B-\Delta B=-i\Delta^{2}(E\times
\bar{E})\\
\end{array}
\right.
\end{equation}
with initial data
\begin{equation}\label{1.2}
E(0,x)=E_{0},\ (n(0,x),n_{t}(0,x))=(n_{0},n_{1}),\
(B(0,x),B_{t}(0,x))=(B_{0},B_{1}),
\end{equation}
where $\alpha\geq 1$ is a constant, $x\in \mathbb{R}^{d}$, $d=2,\
3$. The function $E: \mathbb{R}\oplus \mathbb{R}^{d}\rightarrow
\mathbb{C}^{3}$ is the slowly varying amplitude of the
high-frequency electric field, and the function $n:
\mathbb{R}\oplus \mathbb{R}^{d}\rightarrow \mathbb{R}$ denotes the
fluctuation of the ion-density from its equilibrium, and $B:
\mathbb{R}\oplus \mathbb{R}^{d}\rightarrow \mathbb{R}^{3}$ is the
self-generated magnetic. $\bar{E}$ denotes the conjugate complex
of $E$, and the notation $\times$ appearing in \eqref{1.1} means
the cross product for $\mathbb{R}^{3}$ or $\mathbb{C}^{3}$ valued
vectors. If the space dimension $d=2$, $E$ and $B$ are always
taken as the form $E(t,x)=(E_{1}(t,x),E_{2}(t,x),0)$,
$B(t,x)=(0,0,B_{3}(t,x))$, $x\in \mathbb{R}^{2}$.

Omitting the magnetic field $B$, then the system
\eqref{1.1}-\eqref{1.2} reduces to the standard Zakharov system
(taking $\alpha=1$)
\begin{equation}\label{1.3}
\left\{\!\!
\begin{array}{lc}
iE_{t}+\Delta E=nE, &\\
n_{tt}-\Delta n=\Delta |E|^{2},&\\
E(0,x)=E_{0},\ n(0,x)=n_{0},\ n_{t}(0,x)=n_{1}.
\end{array}
\right.
\end{equation}
This system has been studied by many mathematicians in the past
decades. For the Zakharov system \eqref{1.3}, local existence and
uniqueness of smooth solution $(E,n)\in L^{\infty}(0,T;H^{m}\oplus
H^{m-1})$ with $m\geq 3$ integer were first obtained by C. Sulem
and P. L. Sulem \cite{SS}, in which they also proved the solution
is global in time in one spatial dimension. We also refer to
\cite{GS} for  the results of classical solution in one space
dimensional case. In \cite{AA1}, H. Added and S. Added proved the
smooth solution can be extended globally in time when
$\|E_{0}\|_{L^{2}}$ is small in the case $d=2$. Local
well-posedness in $H^{2}\oplus H^{1}\oplus L^{2}$ was shown by  T.
Ozawa and Y. Tsutsumi in \cite{OT}. J. Bourgain and J. Colliander
\cite{BC} obtained local wellposed results in the energy norm
$(E_{0},n_{0},n_{1})\in H^{1}\oplus L^{2}\oplus \dot{H}^{-1}$ and
showed the solution is global under small assumption on $E_{0}$ in
$d=2,\ 3$. Furthermore, J. Ginibre, Y. Tsutsumi and G. Velo
\cite{GTV} established local well-posedness theory in lower
regularity Sobolev spaces. For more well-posedness results for the
Zakharov system \eqref{1.3}, we refer to \cite{BHHT,CHT,Hol,Pec}
and the references therein.

However, the system \eqref{1.3} ignores the effect of the magnetic
filed which is generated in the laser plasma. In fact, it is
meaningful to consider the self-generated magnetic field in the
Zakharov system from physical viewpoint, e.g. we can study whether
the magnetic field can promote the formation of soliton in three
dimensions or whether it can affect the collapse process of wave
packet in plasma. The magnetic $B$ has  has different expressions
in different plasmas. In a cold plasma, the spontaneous of a
magnetic filed is given by
\begin{equation}\label{1.4}
\Delta B-i\eta \nabla\times \big(\nabla\times (E\times
\bar{E})\big)+\beta B=0,\ \beta\leq 0,\ \eta>0
\end{equation}
while in a hot plasma, the magnetic filed satisfies
\begin{equation}\label{1.5}
\Delta B-i\eta \nabla\times \big(\nabla\times (E\times
\bar{E})\big)-\gamma\frac{\partial}{\partial
t}\int_{\mathbb{R}^{3}}\frac{B(t,y)}{|x-y|^{2}}dy=0,\ \eta,\gamma>
0.
\end{equation}
One can see \cite{KST} for the derivation of the above magnetic
equation. In \cite{Lau}, C. Laurey studied the existence and
uniqueness of the solution for the Zakharov system with the
magnetic given by \eqref{1.4} or \eqref{1.5}. Starting from
Vlasov-Maxwell equations, X. He \cite{He} first derived the
magnetic system \eqref{1.1}, for which describes the pondermotive
force and magnetic field generation effects resulting from the
non-linear interaction between plasma-wave and particles. Hence,
in the present paper, we are devoted to studying the Cauchy
problem of the magnetic Zakharov system \eqref{1.1}-\eqref{1.2}.

To obtain local well-posedness of the the magnetic Zakharov system
\eqref{1.1}-\eqref{1.2}, we use  the energy method together with
communicator estimate to derive a priori bounds and an
approximation argument for the construction of solutions. After
obtaining the uniform bounds for the approximating solutions, we
prove strong convergence of these solutions, then we can get the
well-posedness results. Now we state our main results.

\begin{theorem}\label{thm1.1}
Assume that $s> \frac{d}{2}$, and let $D_{R}(0)$ be the set of
$(E_{0},n_{0},n_{1},B_{0},B_{1})\in H^{s+1}\oplus H^{s}\oplus
(H^{s-1}\cap \dot{H}^{-1})\oplus (H^{s}\cap \dot{H}^{-1})\oplus
(H^{s-2}\cap \dot{H}^{-2})$ such that
$$
\|E_{0}\|_{H^{s+1}}+\|n_{0}\|_{H^{s}}+\|n_{1}\|_{H^{s-1}\cap
\dot{H}^{-1}}+\|B_{0}\|_{H^{s}\cap
\dot{H}^{-1}}+\|B_{1}\|_{H^{s-2}\cap \dot{H}^{-2}}\leq R.
$$
Then for all $R>0$, there exists
$T_{\mathrm{max}}=T_{\mathrm{max}}(R)>0$ such that for all
$(E_{0},n_{0},n_{1},B_{0},B_{1})\in D_{R}(0)$ the magnetic
Zakharov system \eqref{1.1} has a unique solution $(E,n,B)$ with
\begin{equation}\label{1.6}
\begin{array}{cl}
(E,n,B)\in C\big([0,T_{\mathrm{max}});H^{s+1}\oplus H^{s}\oplus
(H^{s }\cap
\dot{H}^{-1})\big),&\\
(E_{t},n_{t},B_{t})\in C\big([0,T_{\mathrm{max}});H^{s-1}\oplus
(H^{s-1}\cap \dot{H}^{-1})\oplus (H^{s-2}\cap \dot{H}^{-2})\big).&
\end{array}
\end{equation}
\end{theorem}

 Note that the above theorem needs the additional condition $n_{1}\in
\dot{H}^{-1}$, $B_{0}\in \dot{H}^{-1}$ and $B_{1}\in
\dot{H}^{-2}$. Since $\mathcal{S}(\mathbb{R}^{2})\not\subset
\dot{H}^{-1}(\mathbb{R}^{2})$ and
$\mathcal{S}(\mathbb{R}^{3})\not\subset
\dot{H}^{-2}(\mathbb{R}^{3})$, the additional assumption on
$n_{1}$, $B_{0}$, $B_{1}$ seems unnatural. In fact, inspired by
\cite{GM1}, this condition can be removed by splitting the initial
data into low frequency part and high frequency part. Namely, we
have the following result concerning the Cauchy problem for the
magnetic Zakharov system.

\begin{theorem}\label{thm1.2}
Assume $s> \frac{d}{2}$, and let $\tilde{D}_{R}(0)$ be the set of
$(E_{0},n_{0},n_{1},B_{0},B_{1})\in H^{s+1}\oplus H^{s}\oplus
H^{s-1}\oplus H^{s}\oplus H^{s-2}$ such that
$$
\|E_{0}\|_{H^{s+1}}+\|n_{0}\|_{H^{s}}+\|n_{1}\|_{H^{s-1}}+\|B_{0}\|_{H^{s}}+\|B_{1}\|_{H^{s-2}}\leq
R.
$$
Then for all $R>0$, there exists
$T_{\mathrm{max}}=T_{\mathrm{max}}(R)>0$ such that for all
$(E_{0},n_{0},n_{1},B_{0},B_{1})\in \tilde{D}_{R}(0)$ the magnetic
Zakharov system \eqref{1.1} has a unique solution $(E,n,B)$ with
\begin{align*}
&(E,n,B)\in C([0,T_{\mathrm{max}});H^{s+1}\oplus H^{s}\oplus H^{s }),\\
&(E_{t},n_{t},B_{t})\in C([0,T_{\mathrm{max}});H^{s-1}\oplus
H^{s-1}\oplus H^{s-2}).
\end{align*}
\end{theorem}

Throughout the paper, the square root of the Laplacian
$(-\Delta)^{\frac{1}{2}}$ will be denoted by $\Lambda$ and
obviously
$$
\mathcal{F}(\Lambda f)=|\xi|\hat{f}.
$$
We denote the inner product of $f$ and $g$ by
$(f,g):=\int_{\mathbb{R}^{d}}f(x)\cdot \bar{g}(x)dx$. We define,
for $s\in \mathbb{R}$ and $1\leq p\leq \infty$, the inhomogeneous
Sobolev space $H^{s,p}(\mathbb{R}^{d})$ or simply $H^{s,p}$ of
tempered distribution $f$ such that
$$
\|f\|_{H^{s,p}}=\|(I-\Delta)^{\frac{s}{2}}f\|_{L^{p}}<\infty,
$$
where $(I-\Delta)^{\frac{s}{2}}f$ is defined by $
(I-\Delta)^{\frac{s}{2}}f=\mathcal{F}^{-1}\Big((1+|\xi|^{2})^{\frac{s}{2}}\hat{f}\Big).
$ If $p=2$, we write $H^{s}$ instead of $H^{s,2}$ for short, and
by Plancherel's theorem $
\|f\|_{H^{s}}=\|(1+|\xi|^{2})^{\frac{s}{2}}\hat{f}\|_{L^{2}}. $
For $s\in \mathbb{R}$, one can define the homogeneous Sobolev
space $\dot{H}^{s,p}$ or $\dot{H}^{s}(=\dot{H}^{s,2})$ by
$$
\|f\|_{\dot{H}^{s,p}}=\|\Lambda^{s}f\|_{L^{p}}<\infty.
$$

This paper is organized as follows. In the next section, we derive
some conserved quantities of the system \eqref{1.1}, and present
the existence result of weak solutions. In Section 3, we introduce
a regularized system for our  magnetic Zakharov system that exists
a unique smooth solution globally. We derive \emph{a priori}
estimates for this regularized system in Section 4 and obtain the
strong convergence property of the approximating solution in
Section 5. Section 6 is concerned with the proof of the main
theorem.

%
%
%
%
%

\section{Conserved quantities and weak solutions}
\setcounter{section}{2}\setcounter{equation}{0}

As we know, conserved laws paly an important role in the analytic
theory(e.g. well-posedness theory and asymptotic behavior) for
nonlinear PDEs of physical origin. For the magnetic Zakharov
system \eqref{1.1}, we have the following conserved results.

\begin{proposition}\label{prop2.1}
For sufficiently regular solutions of the system \eqref{1.1},
there hold two conserved quantities:
\begin{align}
\Phi(t):=&\| E(t)\|_{L^{2}}^{2}=\Phi(0),\label{2.1}\\
\Psi(t):=&\|\nabla\cdot E(t)\|_{L^{2}}^{2}+\alpha\|\nabla\times
E(t)\|_{L^{2}}^{2}+\frac{1}{2}\|n(t)\|_{L^{2}}^{2}+\frac{1}{2}\|\Lambda^{-1}n_{t}(t)\|_{L^{2}}^{2}\nonumber\\
&+
\frac{1}{2}\|\Lambda^{-2}B_{t}(t)\|_{L^{2}}^{2}+\frac{1}{2}\|B(t)\|_{L^{2}}^{2}+\frac{1}{2}\|\Lambda^{-1}B(t)\|_{L^{2}}^{2}\nonumber\\
&+\int_{\mathbb{R}^{d}}n(t)|E(t)|^{2}dx+i\int_{\mathbb{R}^{d}}\big(E(t)\times
\overline{E(t)}\big)\cdot B(t)dx\nonumber\\
&=\Psi(0)\label{2.2}.
\end{align}
\end{proposition}

\begin{proof} Multiplying the first equation of \eqref{1.1} by
$\overline{E}$, then integrating the imaginary part over
$\mathbb{R}^{d}$, and noticing that
\begin{align*}
&2\mathrm{Im}(iE_{t}\cdot \overline{E})=|E|^{2}_{t},\
(\nabla(\nabla\cdot E),E)=-\|\nabla\cdot E\|_{L^{2}}^{2},\ (nE,E)=\|n|E|^{2}\|_{L^{1}},\\
&(\nabla\times(\nabla\times E),E)=\|\nabla\times E\|_{L^{2}}^{2},\
i(E\times B)\cdot \overline{E}=i(\overline{E}\times E)\cdot B
\end{align*}
and $E\times \overline{E}$ is purely imaginary, we then obtain
$$
\frac{1}{2}\frac{d}{dt}\|E(t)\|_{L^{2}}^{2}=0
$$
from which \eqref{2.1} follows.

Now multiplying the first equation of \eqref{1.1} by
$-\overline{E_{t}}$ and integrating the real part, then we have
\begin{align}\label{2.3}
\frac{1}{2}\frac{d}{dt}(\|\nabla\cdot
E\|_{L^{2}}^{2}+\alpha\|\nabla\times
E\|_{L^{2}}^{2})+\frac{1}{2}\int_{\mathbb{R}^{d}}n|E|^{2}_{t}dx-\mathrm{Re}\
i\int_{\mathbb{R}^{d}}(E\times B)\cdot \overline{E_{t}}dx=0.
\end{align}
We take inner product of the second equation of \eqref{1.1} with
$\Lambda^{-2}n_{t}$ and obtain
\begin{align}\label{2.4}
\frac{1}{2}\frac{d}{dt}(\|\Lambda^{-1}n_{t}\|_{L^{2}}^{2}+\|n\|_{L^{2}}^{2})+\int_{\mathbb{R}^{d}}n_{t}|E|^{2}dx=0.
\end{align}
Similarly, if one takes inner product of the third equation of
\eqref{1.1} with $\Lambda^{-4}B_{t}$, then one has
\begin{align}\label{2.5}
\frac{1}{2}\frac{d}{dt}(\|\Lambda^{-2}B_{t}\|_{L^{2}}^{2}+\|B\|_{L^{2}}^{2}+\|\Lambda^{-1}B\|_{L^{2}}^{2})+i\int_{\mathbb{R}^{d}}(E\times
\bar{E})\cdot B_{t}dx=0.
\end{align}
Since $(E\times \bar{E})_{t}=2i \mathrm{Im}(E\times
\overline{E_{t}})$, we then have
\begin{align}
-\mathrm{Re}\ i\int_{\mathbb{R}^{d}}(E\times B)\cdot
\overline{E_{t}}dx&=-\mathrm{Im}\int_{\mathbb{R}^{d}}(E\times
\overline{E_{t}})\cdot B dx\nonumber\\
&=\frac{i}{2} \int_{\mathbb{R}^{d}}(E\times \overline{E})_{t}\cdot
B dx.\label{2.6}
\end{align}
Combining the equalities \eqref{2.3}-\eqref{2.6}, we then get
$$
\frac{d}{dt}\Psi(t)=0
$$
which implies that $\Psi(t)=\Psi(0)$.
\end{proof}

The conserved quantities \eqref{2.1}-\eqref{2.2} are the mail tool
in establishing the global existence of weak solutions for the
system \eqref{1.1}-\eqref{1.2}. Before doing so, we first give the
following elementary lemma.

\begin{lemma}\label{lem2.1} Assume $f(t)$ is a nonnegative
continuous function in $\mathbb{R}^{+}$, and satisfies
$$
f(t)\leq a+bf^{\kappa}(t),\ a,\ b>0,\ \kappa>1.
$$
If $a^{\kappa-1}b<\frac{(\kappa-1)^{\kappa-1}}{\kappa^{\kappa}}$
and $f(0)\leq a$, then $f(t)$ is bounded in $\mathbb{R}^{+}$.
\end{lemma}

\begin{proof}
Let $g(x)=a+bx^{\kappa}-x$, $x\geq 0$. It is easy to see that the
function $g$ has a unique critical point
$x_{0}=\left(\frac{1}{b\kappa}\right)^{\frac{1}{\kappa-1}}$.
Hence, the condition
$a^{\kappa-1}b<\frac{(\kappa-1)^{\kappa-1}}{\kappa^{\kappa}}$
ensures $f(x_{0})<0$ which implies that there exist two points
$x_{1}<x_{2}$ such that $g(x_{1})=g(x_{2})=0$. So if $g(x)\geq 0$,
then either $x\geq x_{2}$ or $0\leq x\leq x_{1}$. Now set
$x=f(t)$, since $a<x_{1}$ and $f(t)$ is continuous, then the
another condition $f(0)\leq a$ ensures $f(t)\leq x_{1}$ for all
$t\geq 0$. Lemma \ref{lem2.1} then follows.
\end{proof}

\begin{lemma}\label{lem2.2}
 Let $(E,n,B)$ be a sufficiently regular solution to the magnetic
Zakharov system \eqref{1.1}-\eqref{1.2} with $\|E_{0}\|_{L^{2}}$
small $(d=2)$ or $\|E_{0}\|_{H^{1}}$ small $(d=3)$, more
precisely,
\begin{align}\label{A2.1}
\begin{split}
&2\|E_{0}\|_{L^{2}}^{2}<\|Q\|_{L^{2}}^{2},\ \mathrm{if}\ d=2,\\
& \|E_{0}\|_{L^{2}}^{2}<\frac{1}{27K^{8}(3)|\Psi(0)|},\ \|\nabla
E_{0}\|_{L^{2}}^{2}\leq |\Psi(0)|,\ \mathrm{if}\ d=3,
\end{split}
\end{align}
where $Q=Q(x)$ is the ground state solution of
\begin{equation*}
\Delta Q-Q+Q^{3}=0,\ x\in \mathbb{R}^{2}.
\end{equation*}
Then we have
\begin{equation}\label{2.7}
\|E\|_{H^{1}}^{2}+\|n\|_{L^{2}}^{2}+\|n_{t}\|_{\dot{H}^{-1}}^{2}+\|B\|_{L^{2}\cap\dot{H}^{-1}}^{2}+\|B_{t}\|_{\dot{H}^{-2}}^{2}\leq
C,
\end{equation}
here the constant $C$ depends on
$\|E_{0}\|_{H^{1}},\|n_{0}\|_{L^{2}},
\|n_{1}\|_{\dot{H}^{-1}},\|B_{0}\|_{L^{2}\cap\dot{H}^{-1}},\|B_{1}\|_{\dot{H}^{-2}}$.
\end{lemma}

\begin{proof} By Cauchy-Schwarz inequality, we have
\begin{equation}\label{2.8}
\left|\int_{\mathbb{R}^{d}}n|E|^{2}dx\right|\leq \epsilon
\|n\|_{L^{2}}^{2}+\frac{1}{4\epsilon}\|E\|_{L^{4}}^{4}\leq
\epsilon
\|n\|_{L^{2}}^{2}+\frac{1}{4\epsilon}K^{4}(d)\|E\|_{L^{2}}^{4-d}\|\nabla
E\|^{d}_{L^{2}}
\end{equation}
for all $0<\epsilon<\frac{1}{2}$, where we have used the following
Sobolev best constant inequality (see \cite{Wei})
$$
\|f\|_{L^{4}(\mathbb{R}^{d})}^{4}\leq
K^{4}(d)\|f\|_{L^{2}(\mathbb{R}^{d})}^{4-d}\|\nabla
f\|^{d}_{L^{2}(\mathbb{R}^{d})},
$$
and $K^{4}(d)=\frac{2}{\|Q\|_{L^{2}}^{2}}$, $Q$ is the ground
state solution of
$$
\frac{d}{2}\Delta Q-(2-\frac{d}{2})Q+Q^{3}=0.
$$
Similarly, we can obtain
\begin{equation}\label{2.9}
\left|i\int_{\mathbb{R}^{d}}\big(E(t)\times
\overline{E(t)}\big)\cdot B(t)dx\right|\leq \epsilon
\|B\|_{L^{2}}^{2}+\frac{1}{4\epsilon}K^{4}(d)\|E\|_{L^{2}}^{4-d}\|\nabla
E\|^{d}_{L^{2}}
\end{equation}
for all $0<\epsilon<\frac{1}{2}$.

Due to the fact $\|\nabla\cdot E(t)\|_{L^{2}}^{2}+\|\nabla\times
E(t)\|_{L^{2}}^{2}=\|\nabla E\|_{L^{2}}^{2}$, we deduce from
\eqref{2.8}, \eqref{2.9} and Proposition \ref{prop2.1} that if
$d=2$
\begin{align*}
&\|\nabla E\|_{L^{2}}^{2}+(\frac{1}{2}-\epsilon)\|n\|_{L^{2}}^{2}
+(\frac{1}{2}-\epsilon)\|B\|_{L^{2}}^{2}+\frac{1}{2}\|B\|_{\dot{H}^{-1}}^{2}+\frac{1}{2}\|n_{t}\|_{\dot{H}^{-1}}^{2}+\frac{1}{2}\|B_{t}\|_{\dot{H}^{-2}}^{2}\\
&\quad\quad\leq
|\Psi(0)|+\frac{1}{2\epsilon}K^{4}(2)\|E_{0}\|_{L^{2}}^{2}\|\nabla
E\|^{2}_{L^{2}}\\
&\quad\quad=
|\Psi(0)|+\frac{1}{\epsilon}\frac{\|E_{0}\|_{L^{2}}^{2}}{\|Q\|_{L^{2}}^{2}}\|\nabla
E\|^{2}_{L^{2}}.
\end{align*}
If $2\|E_{0}\|_{L^{2}}^{2}<\|Q\|_{L^{2}}^{2}$, then we can choose
$\epsilon$ very close to $\frac{1}{2}$ such that
$\|E_{0}\|_{L^{2}}^{2}<\epsilon\|Q\|_{L^{2}}^{2}$,  thus
\eqref{2.7} follows in the case $d=2$.

If $d=3$, we can obtain($\epsilon=\frac{1}{4}$)
\begin{align}
&\|\nabla E\|_{L^{2}}^{2}+\frac{1}{4}\|n\|_{L^{2}}^{2}
+\frac{1}{4}\|B\|_{L^{2}}^{2}+\frac{1}{2}\|B\|_{\dot{H}^{-1}}^{2}+\frac{1}{2}\|n_{t}\|_{\dot{H}^{-1}}^{2}+\frac{1}{2}\|B_{t}\|_{\dot{H}^{-2}}^{2}\nonumber\\
&\quad\quad\leq |\Psi(0)|+2K^{4}(3)\|E_{0}\|_{L^{2}}\|\nabla
E\|^{3}_{L^{2}}.\label{2.10}
\end{align}
If we take $f(t)=\|\nabla E\|_{L^{2}}^{2}$, $a= |\Psi(0)|$,
$b=2K^{4}(3)\|E_{0}\|_{L^{2}}$, $\kappa=\frac{3}{2}$, then Lemma
\ref{lem2.1} gives $\|\nabla E(t)\|_{L^{2}}^{2}\leq C$ for all
$t\geq 0$. Hence, \eqref{2.7} is obtained from \eqref{2.10}.
\end{proof}

So an immediate application of the conservation laws
\eqref{2.1}-\eqref{2.2} is to establish the existence of weak
solutions to the magnetic Zakharov system \eqref{1.1}.

\begin{theorem}\label{thm2.1}
If $E_{0}\in H^{1}$, $(n_{0},n_{1})\in L^{2}\oplus \dot{H}^{-1}$,
$(B_{0},B_{1})\in (L^{2}\cap \dot{H}^{-1})\oplus \dot{H}^{-2}$,
and the initial data satisfying \eqref{A2.1}, then there exists a
weak solution $(E,n,B)$ for the system \eqref{1.1}-\eqref{1.2} in
the distributional sense such that
$$
E\in L^{\infty}(\mathbb{R}^{+};H^{1}),\ (n,n_{t})\in
L^{\infty}(\mathbb{R}^{+};L^{2}\oplus \dot{H}^{-1}),\ (B,B_{t})\in
L^{\infty}(\mathbb{R}^{+};(L^{2}\cap
\dot{H}^{-1})\oplus\dot{H}^{-2}).
$$
\end{theorem}

Using the prior estimate \eqref{2.7}, Theorem \ref{thm2.1} can be
proved by applying Galerkin method and compactness argument, since
this procedure is standard, the proof of Theorem \ref{thm2.1} is
omitted here.

%
%
%
%
%

\section{Regularization for the original system}
\setcounter{section}{3}\setcounter{equation}{0}

In this section, we introduce a regularized system for our
original system \eqref{1.1}. Now consider the following
system($0<\epsilon<1$)
\begin{subequations}
\begin{align}
&
iE_{t}^{\epsilon}+i\epsilon^{2}\Delta^{2}E_{t}^{\epsilon}+\nabla(\nabla\cdot
E^{\epsilon})-\alpha\nabla\times(\nabla\times
E^{\epsilon})-n^{\epsilon}E^{\epsilon}+i E^{\epsilon}\times B^{\epsilon}=0, \label{3.1a}\\
&n^{\epsilon}_{tt}-\Delta n^{\epsilon}=\Delta |E^{\epsilon}|^{2},\label{3.1b}\\
&B^{\epsilon}_{tt}+\Delta ^{2}B^{\epsilon}-\Delta
B^{\epsilon}=-i\Delta^{2}(E^{\epsilon}\times
\overline{E^{\epsilon}})\label{3.1c}
 \end{align}
with smooth initial data
\begin{equation}\label{3.1d}
E^{\epsilon}(0)=E^{\epsilon}_{0},\
n^{\epsilon}(0)=n^{\epsilon}_{0},\
n^{\epsilon}_{t}(0)=n^{\epsilon}_{1},\
B^{\epsilon}(0)=B^{\epsilon}_{0},\
B^{\epsilon}_{t}(0)=B^{\epsilon}_{1}.
\end{equation}
\end{subequations}

With the same argument as Proposition \ref{prop2.1}, we can obtain
some conservation results for this regularized system.

\begin{proposition}\label{prop3.1}
Assume $(E^{\epsilon},n^{\epsilon},B^{\epsilon})$ is a sufficient
regular solution for the system \eqref{3.1a}-\eqref{3.1d}, then we
have
\begin{align}
\Phi^{\epsilon}(t):=&\| E^{\epsilon}(t)\|_{L^{2}}^{2}+ \epsilon^{2}\| \Delta E^{\epsilon}(t)\|_{L^{2}}^{2}=\Phi^{\epsilon}(0),\label{3.2}\\
\Psi^{\epsilon}(t):=&\|\nabla\cdot
E^{\epsilon}(t)\|_{L^{2}}^{2}+\alpha\|\nabla\times
E^{\epsilon}(t)\|_{L^{2}}^{2}+\frac{1}{2}\|n^{\epsilon}(t)\|_{L^{2}}^{2}+\frac{1}{2}\|\Lambda^{-1}n^{\epsilon}_{t}(t)\|_{L^{2}}^{2}\nonumber\\
&+
\frac{1}{2}\|\Lambda^{-2}B^{\epsilon}_{t}(t)\|_{L^{2}}^{2}+\frac{1}{2}\|B^{\epsilon}(t)\|_{L^{2}}^{2}
+\frac{1}{2}\|\Lambda^{-1}B^{\epsilon}(t)\|_{L^{2}}^{2}\nonumber\\
&+\int_{\mathbb{R}^{d}}n^{\epsilon}(t)|E^{\epsilon}(t)|^{2}dx+i\int_{\mathbb{R}^{d}}\big(E^{\epsilon}(t)\times
\overline{E^{\epsilon}(t)}\big)\cdot B^{\epsilon}(t)dx\nonumber\\
&=\Psi^{\epsilon}(0)\label{3.3}.
\end{align}
\end{proposition}

Let $\mathcal{L}=(I+\epsilon^{2}\Delta ^{2})^{-1}$, and let
$\mathcal{A}$ be the linear operator defined by
$$
\mathcal{A}E=-\nabla(\nabla\cdot E)+\alpha \nabla\times
(\nabla\times E),
$$
since the operator $\mathcal{L}\mathcal{A}$ is self-adjoint, then
the linear equation
$$
iE_{t}=\mathcal{L}\mathcal{A}E,\ E(0)=E_{0}
$$
generates a unitary group $U(t)$ in $H^{r}(\mathbb{R}^{d})$.
Therefore, we can transform the regularized system
\eqref{3.1a}-\eqref{3.1d} into the following integral equation
\begin{align}\label{3.4}
E^{\epsilon}(t)=U(t)E^{\epsilon}_{0}+\int_{0}^{t}U(t-\tau)f(E^{\epsilon}(\tau))d\tau,
\end{align}
where
$f(E^{\epsilon}(t))=-i\mathcal{L}(n^{\epsilon}E^{\epsilon})-\mathcal{L}(
E^{\epsilon}\times B^{\epsilon})$, and
$n^{\epsilon}=n^{\epsilon}(E^{\epsilon})$,
$B^{\epsilon}=B^{\epsilon}(E^{\epsilon})$ is the solution of
equation \eqref{3.1b}, \eqref{3.1c} respectively. Indeed, we can
express $n^{\epsilon}(E^{\epsilon})$, $B^{\epsilon}(E^{\epsilon})$
by
\begin{align}
n^{\epsilon}(E^{\epsilon})&=\cos((-\Delta)^{1/2}t)n_{0}^{\epsilon}
+\frac{\sin((-\Delta)^{1/2}t)}{(-\Delta)^{1/2}}n_{1}^{\epsilon}\nonumber\\
&\quad+\int_{0}^{t}\frac{\sin((-\Delta)^{1/2}(t-\tau))}{(-\Delta)^{1/2}}\Delta
|E^{\epsilon}(\tau)|^{2}d\tau,\label{3.5}
\end{align}
and
\begin{align}
B^{\epsilon}(E^{\epsilon})&=\cos((-\Delta)^{1/2}\langle\nabla\rangle
t)B_{0}^{\epsilon}
+\frac{\sin((-\Delta)^{1/2}\langle\nabla\rangle t)}{(-\Delta)^{1/2}\langle\nabla\rangle}B_{1}^{\epsilon}\nonumber\\
&\quad-i\int_{0}^{t}\frac{\sin((-\Delta)^{1/2}\langle\nabla\rangle(t-\tau))}{(-\Delta)^{1/2}\langle\nabla\rangle}\Delta^{2}
(E^{\epsilon}(\tau)\times
\overline{E^{\epsilon}(\tau)})d\tau,\label{3.6}
\end{align}
here $\langle\nabla\rangle$ is the Fourier multiplier with symbol
$\langle\xi\rangle=(1+|\xi|^{2})^{1/2}$.

The main result in this section is the following global existence
of smooth solution for the regularized system
\eqref{3.1a}-\eqref{3.1d}.

\begin{theorem}\label{thm3.1}
Given $\epsilon\in (0,1)$, and suppose that $E^{\epsilon}_{0}\in
H^{r+1}$,\ $(n^{\epsilon}_{0},n^{\epsilon}_{1})\in H^{r}\oplus
(H^{r-1}\cap \dot{H}^{-1})$,
$(B^{\epsilon}_{0},B^{\epsilon}_{1})\in
(H^{r-1}\cap\dot{H}^{-1})\oplus (H^{r-3}\cap \dot{H}^{-2}))$, $r$
is large enough$($e.g. $r>\frac{d}{2}+s+100)$, then there exists a
unique smooth solution $(E^{\epsilon},n^{\epsilon},B^{\epsilon})$
for the regularized system \eqref{3.1a}-\eqref{3.1d} such that
\begin{align*}
(E^{\epsilon},n^{\epsilon},B^{\epsilon})\in
C(\mathbb{R}^{+};H^{r+1}\oplus H^{r}\oplus H^{r-1}).
\end{align*}
\end{theorem}

In order to prove Theorem \ref{thm3.1}, we first state the
following calculus inequality which will be used many times in
this paper.

\begin{lemma}\label{lem1}
Assume that $s>0$ and $p\in (1,+\infty)$. If $f,g\in
\mathcal{S}(\mathbb{R}^{d})$ , the Schwartz class, then
\begin{equation}\label{3.7}
\|\Lambda^{s}(fg)\|_{L^{p}}\leq
C(\|f\|_{L^{p_{1}}}\|g\|_{\dot{H}^{s,p_{2}}}+\|f\|_{\dot{H}^{s,p_{3}}}\|g\|_{L^{p_{4}}})
\end{equation}
and
\begin{equation}\label{3.8}
\|\Lambda^{s}(fg)-f(\Lambda^{s}g)\|_{L^{p}}\leq C\|\nabla
f\|_{L^{p_{1}}}\|g\|_{\dot{H}^{s-1,p_{2}}}+\|f\|_{\dot{H}^{s,p_{3}}}\|g\|_{L^{p_{4}}}
\end{equation}
 with $p_{2},p_{3}\in (1,+\infty)$ such that
\begin{equation*}
\frac{1}{p}=\frac{1}{p_{1}}+\frac{1}{p_{2}}=\frac{1}{p_{3}}+\frac{1}{p_{4}}.
\end{equation*}
\end{lemma}

For a proof of this lemma, we refer to \cite{CM,K}.

{\it Proof of Theorem \ref{thm3.1}.}  By contraction argument, we
first show that equation \eqref{3.4} has a unique solution
locally, then we extend this solution globally in time based on
some uniform estimates.

Denote
$M=\|E_{0}^{\epsilon}\|_{H^{r+1}}+\|n_{0}^{\epsilon}\|_{H^{r}}+\|n_{1}^{\epsilon}\|_{H^{r-1}}+\|B_{0}^{\epsilon}\|_{H^{r-1}}
+\|B_{1}^{\epsilon}\|_{H^{r-3}}$. Let $0<T\leq1$ be determined
later, and set $X=C([0,T];H^{r+1})$, $X_{M}=\{E\in X;
\|E\|_{X}\leq 2M\}$. Now we define the map $\mathcal{T}$ acting on
$X$ by
\begin{align}\label{3.9}
\mathcal{T}(E^{\epsilon})=U(t)E^{\epsilon}_{0}+\int_{0}^{t}U(t-\tau)f(E^{\epsilon}(\tau))d\tau,\
\forall\ E^{\epsilon}\in X.
\end{align}
 Our aim is to show that
$\mathcal{T}$ has unique fixed point on $X$ if $T$ is small
enough.

Given $E^{\epsilon}\in X$, we obtain from \eqref{3.5} that
\begin{align}\label{3.10}
\|n^{\epsilon}(t)\|_{H^{r}}\leq
(1+T)(\|n_{0}^{\epsilon}\|_{H^{r}}+\|n_{1}^{\epsilon}\|_{H^{r-1}})+T\|E^{\epsilon}\|_{X}^{2},
\forall\ t\in [0,T].
\end{align}
Similarly, from \eqref{3.6} we have
\begin{align}\label{3.11}
\|B^{\epsilon}(t)\|_{H^{r-1}}\leq
(1+T)(\|B_{0}^{\epsilon}\|_{H^{r-1}}+\|n_{1}^{\epsilon}\|_{H^{r-1}})+T\|E^{\epsilon}\|_{X}^{2},
\forall\ t\in [0,T].
\end{align}
Note also that $\mathcal{L}$ is a bounded linear operator from
$H^{k-4}$ to $H^{k}$, namely, there exists $K>0$ not depending on
$k$ such that
$$
\|\mathcal{L}f\|_{H^{k}}\leq K\|f\|_{H^{k-4}},\ \forall\ f\in
H^{k-4}.
$$
This fact together with \eqref{3.10}-\eqref{3.11} yield
\begin{align*}
\|f(E^{\epsilon})\|_{X}&\leq
K(\|n^{\epsilon}(E^{\epsilon})E^{\epsilon}\|_{C([0,T];H^{r-3})}+\|E^{\epsilon}\times
B^{\epsilon}(E^{\epsilon})\|_{C([0,T];H^{r-3})})\\
& \leq K\sup\limits_{t\in
[0,T]}(\|n^{\epsilon}(E^{\epsilon})\|_{H^{r}}\|E^{\epsilon}\|_{H^{r+1}}+\|B^{\epsilon}(E^{\epsilon})\|_{H^{r-1}}\|E^{\epsilon}\|_{H^{r+1}})\\
& \leq 4KM^{2}(1+4M),\ \mathrm{if}\ E^{\epsilon}\in X_{M},
\end{align*}
then we have
\begin{align*}
\|\mathcal{T}(E^{\epsilon})\|_{X}&\leq
\|E^{\epsilon}_{0}\|_{H^{r+1}}+\int_{0}^{T}\|f(E^{\epsilon})\|_{X}d\tau\\
&\leq M+4KM^{2}(1+4M)T.
\end{align*}
Hence we see that if $T\leq \min\{\frac{1}{4KM(1+4M)},1\}=:T_{1}$,
then $\mathcal{T}$ maps $X_{M}$ into itself.

 From the expression \eqref{3.5} and
\eqref{3.6}, we also have
\begin{align*}
\|n(E^{\epsilon}_{1})(t)-n(E^{\epsilon}_{2})(t)\|_{H^{r}}&\leq
T(\|E^{\epsilon}_{1}\|_{X}+\|E^{\epsilon}_{2}\|_{X})\|E^{\epsilon}_{1}-E^{\epsilon}_{2}\|_{X},\
\forall \ t\in [0,T],\\
\|B(E^{\epsilon}_{1})(t)-B(E^{\epsilon}_{2})(t)\|_{H^{r-1}}&\leq
T(\|E^{\epsilon}_{1}\|_{X}+\|E^{\epsilon}_{2}\|_{X})\|E^{\epsilon}_{1}-E^{\epsilon}_{2}\|_{X},\
\forall \ t\in [0,T].
\end{align*}
Using the above two estimates and \eqref{3.10}-\eqref{3.11}, we
know
\begin{align*}
&\|f(E^{\epsilon}_{1})-f(E^{\epsilon}_{2})\|_{X}\\
&\quad \leq K\sup\limits_{t\in
[0,T]}(\|n^{\epsilon}(E^{\epsilon}_{1})E^{\epsilon}_{1}-n^{\epsilon}(E^{\epsilon}_{2})E^{\epsilon}_{2}\|_{H^{r-3}}
+\|E^{\epsilon}_{1}\times
B^{\epsilon}(E^{\epsilon}_{1})-E^{\epsilon}_{2}\times
B^{\epsilon}(E^{\epsilon}_{2})\|_{H^{r-3}})\\
&\quad \leq K\sup\limits_{t\in
[0,T]}(\|n^{\epsilon}(E^{\epsilon}_{1})\|_{H^{r}}\|E^{\epsilon}_{1}-E^{\epsilon}_{2}\|_{H^{r+1}}
+\|n^{\epsilon}(E^{\epsilon}_{1})-n^{\epsilon}(E^{\epsilon}_{2})\|_{H^{r}}\|E^{\epsilon}_{2}\|_{H^{r+1}}\\
&\quad \quad
+\|B^{\epsilon}(E^{\epsilon}_{1})\|_{H^{r-1}}\|E^{\epsilon}_{1}-E^{\epsilon}_{2}\|_{H^{r+1}}
+\|B^{\epsilon}(E^{\epsilon}_{1})-B^{\epsilon}(E^{\epsilon}_{2})\|_{H^{r-1}}\|E^{\epsilon}_{2}\|_{H^{r+1}})\\
&\quad \leq
K(24M^{2}+2M)\|E^{\epsilon}_{1}-E^{\epsilon}_{2}\|_{X}, \
\mathrm{if}\ E^{\epsilon}_{1}, E^{\epsilon}_{2}\in X_{M}.
\end{align*}
Therefore, there holds
$$
\|\mathcal{T}(E^{\epsilon}_{1})-\mathcal{T}(E^{\epsilon}_{2})\|_{X}\leq
TK(24M^{2}+2M)\|E^{\epsilon}_{1}-E^{\epsilon}_{2}\|_{X}.
$$
If we choose $T\leq \min\{\frac{1}{2K(24M^{2}+2M)},T_{1}\}$, then
$\mathcal{T}$ is a contraction map on $X_{M}$. So by fixed point
theorem, we know that the equation \eqref{3.9} has a unique
solution $E^{\epsilon}\in C([0,T];H^{r+1})$, and by
\eqref{3.5}-\eqref{3.6}, we know $n^{\epsilon}\in C([0,T];H^{r})$,
$B^{\epsilon}\in C([0,T];H^{r-1})$.

Assume now $T_{\mathrm{max}}$ is the maximal existence time of the
solution $(E^{\epsilon},n^{\epsilon},B^{\epsilon})$ for the
regularized system \eqref{3.1a}-\eqref{3.1d}, hence in order to
complete the proof of Theorem \ref{thm3.1}, we have to show that
$T_{\mathrm{max}}=\infty$. To this end, it is sufficient to prove
that the quantities $\|E^{\epsilon}\|_{H^{r+1}}$,
$\|n^{\epsilon}\|_{H^{r}}$, $\|n_{t}^{\epsilon}\|_{H^{r-1}}$,
$\|B^{\epsilon}\|_{H^{r-1}}$, $\|B_{t}^{\epsilon}\|_{H^{r-3}}$ are
bounded on the time interval $[0,T_{\mathrm{max}})$.

It follows from \eqref{3.2} that
\begin{equation}\label{3.12}
\|E^{\epsilon}\|_{H^{2}}\leq C(\Rightarrow
\|E^{\epsilon}\|_{L^{\infty}}\leq C),
\end{equation}
 which in turn
gives that(by \eqref{3.5} and \eqref{3.6})
\begin{equation}\label{3.13}
\|n^{\epsilon}\|_{H^{1}}\leq C,\ \|B^{\epsilon}\|_{L^{2}}\leq C.
\end{equation}
Besides, one also deduces from \eqref{3.3} that
\begin{equation}\label{3.14}
\|n^{\epsilon}_{t}\|_{\dot{H}^{-1}}+\|B^{\epsilon}_{t}\|_{\dot{H}^{-2}}\leq
C.
\end{equation}
We emphasize that the constant $C$ in the above estimates depends
on the parameter $\epsilon$.

Now multiplying \eqref{3.1a} by
$\Lambda^{2r-2}\overline{E^{\epsilon}}$, and integrating the
imaginary part, since
$$
\mathrm{Im}\big(\nabla(\nabla\cdot
E^{\epsilon})-\alpha\nabla\times(\nabla\times
E^{\epsilon}),\Lambda^{2r-2}E^{\epsilon}\big)=0,
$$
then we obtain from \eqref{3.7}, \eqref{3.12}-\eqref{3.13} that
\begin{align}\label{3.15}
\begin{split}
&\frac{1}{2}\frac{d}{dt}(\|\Lambda^{r-1}E^{\epsilon}\|_{L^{2}}^{2}+\epsilon^{2}\|\Lambda^{r+1}E^{\epsilon}\|_{L^{2}}^{2})\\
&\quad=
\int_{\mathbb{R}^{d}}\Lambda^{r-1}(n^{\epsilon}E^{\epsilon})\Lambda^{r-1}\overline{E^{\epsilon}}dx-i\int_{\mathbb{R}^{d}}\Lambda^{r-1}(E^{\epsilon}\times
B^{\epsilon})\Lambda^{r-1}\overline{E^{\epsilon}}dx\\
&\quad \leq
C(\|n^{\epsilon}\|_{H^{r-1}}\|E^{\epsilon}\|_{L^{\infty}}+\|n^{\epsilon}\|_{L^{4}}\|E^{\epsilon}\|_{H^{r-1,4}})\|\Lambda^{r-1}\overline{E^{\epsilon}}\|_{L^{2}}\\
&\quad \quad +
C(\|B^{\epsilon}\|_{H^{r-1}}\|E^{\epsilon}\|_{L^{\infty}}+\|B^{\epsilon}\|_{L^{p}}\|E^{\epsilon}\|_{H^{r-1,q}})\|\Lambda^{r-1}\overline{E^{\epsilon}}\|_{L^{2}}\\
 &\quad \leq
C(\|n^{\epsilon}\|_{H^{r}}+\|B^{\epsilon}\|_{H^{r-1}}+\|E^{\epsilon}\|_{H^{r+1}})\|E^{\epsilon}\|_{H^{r+1}}\\
&\quad\leq
C(\|n^{\epsilon}\|_{H^{r}}^{2}+\|B^{\epsilon}\|_{H^{r-1}}^{2}+\|E^{\epsilon}\|_{H^{r+1}}^{2}),
\end{split}
\end{align}
where we have used  the following
inequality($\frac{1}{p}+\frac{1}{q}=\frac{1}{2}$)
\begin{align*}
\overline{\lim\limits_{p\rightarrow
2^{+}}}\|B^{\epsilon}\|_{L^{p}}\|E^{\epsilon}\|_{H^{r-1,q}}\leq
C\|E^{\epsilon}\|_{H^{r+1}}\overline{\lim\limits_{p\rightarrow
2^{+}}}\|B^{\epsilon}\|_{L^{p}}\leq C\|E^{\epsilon}\|_{H^{r+1}}
\end{align*}
due to the fact $\|B^{\epsilon}\|_{L^{2}}\leq C$.

We then multiply \eqref{3.1b} by $\Lambda^{2r-2}n^{\epsilon}_{t}$
and obtain
\begin{align}\label{3.16}
\begin{split}
&\frac{1}{2}\frac{d}{dt}(\|\Lambda^{r-1}n^{\epsilon}_{t}\|_{L^{2}}^{2}+\|\Lambda^{r}n^{\epsilon}\|_{L^{2}}^{2})\\
&\quad=-\int_{\mathbb{R}^{d}}\Lambda^{r+1}(E^{\epsilon}\cdot\overline{E^{\epsilon}})\Lambda^{r-1}n^{\epsilon}_{t}dx\\
&\quad \leq
C(\|E^{\epsilon}\|_{H^{r+1}}\|\overline{E^{\epsilon}}\|_{L^{\infty}}
+\|E^{\epsilon}\|_{L^{\infty}}\|\overline{E^{\epsilon}}\|_{H^{r+1}})\|\Lambda^{r-1}n^{\epsilon}_{t}\|_{L^{2}}\\
&\quad\leq
C(\|n^{\epsilon}_{t}\|_{H^{r-1}}^{2}+\|E^{\epsilon}\|_{H^{r+1}}^{2}).
\end{split}
\end{align}
Similarly, by taking inner product of \eqref{3.1c} with
$\Lambda^{2r-6}B^{\epsilon}_{t}$, then  the same argument as above
leads to
\begin{align}\label{3.17}
\begin{split}
&\frac{1}{2}\frac{d}{dt}(\|\Lambda^{r-3}B^{\epsilon}_{t}\|_{L^{2}}^{2}+\|\Lambda^{r-1}B^{\epsilon}\|_{L^{2}}^{2}+\|\Lambda^{r-2}B^{\epsilon}\|_{L^{2}}^{2})\\
&\quad=\int_{\mathbb{R}^{d}}\Lambda^{r+1}(E^{\epsilon}\times\overline{E^{\epsilon}})\Lambda^{r-3}n^{\epsilon}_{t}dx\\
&\quad\leq
C(\|B^{\epsilon}_{t}\|_{H^{r-3}}^{2}+\|E^{\epsilon}\|_{H^{r+1}}^{2}).
\end{split}
\end{align}

Now summing the estimates \eqref{3.15}-\eqref{3.17}, and
integrating the result, we can obtain
\begin{align*}
&\|\Lambda^{r-1}E^{\epsilon}\|_{L^{2}}^{2}+\epsilon^{2}\|\Lambda^{r+1}E^{\epsilon}\|_{L^{2}}^{2}
+\|\Lambda^{r-1}n^{\epsilon}_{t}\|_{L^{2}}^{2}+\|\Lambda^{r}n^{\epsilon}\|_{L^{2}}^{2}\\
&\quad\quad\quad\quad+\|\Lambda^{r-3}B^{\epsilon}_{t}\|_{L^{2}}^{2}+\|\Lambda^{r-1}B^{\epsilon}\|_{L^{2}}^{2}+\|\Lambda^{r-2}B^{\epsilon}\|_{L^{2}}^{2}\\
&\quad\quad\leq
C+C\int_{0}^{t}(\|n^{\epsilon}\|_{H^{r}}^{2}+\|B^{\epsilon}\|_{H^{r-1}}^{2}+\|E^{\epsilon}\|_{H^{r+1}}^{2}
+\|n^{\epsilon}_{t}\|_{H^{r-1}}^{2}+\|B^{\epsilon}_{t}\|_{H^{r-3}}^{2})d\tau,
\end{align*}
this inequality together with the fact
$\|E^{\epsilon}\|_{L^{2}}^{2}\leq C$ and \eqref{3.13}-\eqref{3.14}
yield
\begin{align*}
&\|E^{\epsilon}\|_{H^{r+1}}^{2}+\|n^{\epsilon}_{t}\|_{H^{r-1}}^{2}+\|n^{\epsilon}\|_{H^{r}}^{2}+\|B^{\epsilon}_{t}\|_{H^{r-3}}^{2}+\|B^{\epsilon}\|_{H^{r-1}}^{2}\\
&\quad\quad\leq
C+C\int_{0}^{t}(\|E^{\epsilon}\|_{H^{r+1}}^{2}+\|n^{\epsilon}\|_{H^{r}}^{2}+\|n^{\epsilon}_{t}\|_{H^{r-1}}^{2}+\|B^{\epsilon}\|_{H^{r-1}}^{2}
+\|B^{\epsilon}_{t}\|_{H^{r-3}}^{2})d\tau,
\end{align*}
hence by Gronwall's inequality, we get
\begin{equation}\label{3.18}
\|E^{\epsilon}\|_{H^{r+1}}^{2}+\|n^{\epsilon}_{t}\|_{H^{r-1}}^{2}
+\|n^{\epsilon}\|_{H^{r}}^{2}+\|B^{\epsilon}_{t}\|_{H^{r-3}}^{2}+\|B^{\epsilon}\|_{H^{r-1}}^{2}\leq
C,\ \forall\ t\in [0,T_{\mathrm{max}}).
\end{equation}
The estimate \eqref{3.18} implies that the solution
$(E^{\epsilon},n^{\epsilon},B^{\epsilon})$ can be extended to the
interval $[0,T_{\mathrm{max}}+\delta]$, which contradicts the
maximality, hence $T_{\mathrm{max}}=\infty$. Therefore, the
solution for regularized system \eqref{3.1a}-\eqref{3.1d} exists
globally in time, and the proof of Theorem \ref{thm3.1} is
complete. \qed

%
%
%
%
%

\section{\textit{A prior} estimates}
\setcounter{section}{4}\setcounter{equation}{0}

We will approximate the solution of the magnetic Zakharov system
\eqref{1.1}-\eqref{1.2} by smooth solutions for the regularized
system given in Section 3. Hence, in order to get strong or weak
limit of these smooth solutions in the $L^{\infty}_{t}H^{s}_{x}$
topology, one must demonstrate the approximating solutions are
uniformly bounded in this energy norm. Therefore, we are devoted
to establishing \textit{a prior} estimates for the system
\eqref{3.1a}-\eqref{3.1d} in this section.

\begin{proposition}\label{prop4.1}
Let $s>\frac{d}{2}$, $E_{0}\in H^{s+1}$, $(n_{0},n_{1})\in
H^{s}\oplus (H^{s-1}\cap \dot{H}^{-1})$, $(B_{0},B_{1})\in
(H^{s}\cap \dot{H}^{-1})\oplus (H^{s-2}\cap \dot{H}^{-2})$. Assume
the sequence
$\{(E_{0}^{\epsilon},n_{0}^{\epsilon},n_{1}^{\epsilon},B_{0}^{\epsilon},B_{1}^{\epsilon})\}$
satisfying $E^{\epsilon}_{0}\in H^{r+1}$,\
$(n^{\epsilon}_{0},n^{\epsilon}_{1})\in H^{r}\oplus (H^{r-1}\cap
\dot{H}^{-1})$, $(B^{\epsilon}_{0},B^{\epsilon}_{1})\in
(H^{r-1}\cap\dot{H}^{-1})\oplus (H^{r-3}\cap \dot{H}^{-2})$ with
$r$ large enough, and
\begin{equation}\label{4.1}
\begin{array}{cl}
E_{0}^{\epsilon}\rightarrow E_{0}\ \mathrm{in}\ H^{s+1},&\\
n_{0}^{\epsilon}\rightarrow n_{0}\ \mathrm{in}\ H^{s},\
n_{1}^{\epsilon}\rightarrow n_{1}\ \mathrm{in}\
H^{s-1}\cap\dot{H}^{-1},&\\
B_{0}^{\epsilon}\rightarrow B_{0}\ \mathrm{in}\
H^{s}\cap\dot{H}^{-1},\ B_{1}^{\epsilon}\rightarrow B_{1}\
\mathrm{in}\ H^{s-2}\cap\dot{H}^{-2}.
\end{array}
\end{equation}
If $(E^{\epsilon},n^{\epsilon},B^{\epsilon})$ is the smooth
solution of regularized system \eqref{3.1a}-\eqref{3.1c} with the
initial data
$(E_{0}^{\epsilon},n_{0}^{\epsilon},n_{1}^{\epsilon},B_{0}^{\epsilon},B_{1}^{\epsilon})$,
then there exist $T>0$ and $C>0$ such that
\begin{align}
&\|E^{\epsilon}\|_{C([0,T];H^{s+1})}+\|n^{\epsilon}\|_{C([0,T];H^{s})}+\|B^{\epsilon}\|_{C([0,T];H^{s}\cap
\dot{H}^{-1})}\leq C,\label{4.2}\\
&
\|E^{\epsilon}_{t}\|_{C([0,T];H^{s-1})}+\|n^{\epsilon}_{t}\|_{C([0,T];H^{s-1}\cap
\dot{H}^{-1})}+\|B^{\epsilon}_{t}\|_{C([0,T];H^{s-2}\cap
\dot{H}^{-2})}\leq C, \label{A4.1}
\end{align}
where $T$ and $C$ are dependent of the norm of
$(E_{0},n_{0},n_{1},B_{0},B_{1})$, but independent of $\epsilon$.
\end{proposition}

\begin{proof}
Note that the bound for the left hand side of \eqref{3.18} depends
on $\epsilon$, hence, one can not use the same argument that leads
to \eqref{3.18} to obtain the estimate \eqref{4.2}-\eqref{A4.1}.
In order to derive independent of $\epsilon$, we first write
\eqref{3.1a} in the following form
\begin{equation}\label{4.3}
iE_{t}^{\epsilon}+\mathcal{L}\nabla(\nabla\cdot
E^{\epsilon})-\alpha\mathcal{L}\nabla\times(\nabla\times
E^{\epsilon})-\mathcal{L}(n^{\epsilon}E^{\epsilon})+i\mathcal{L}(
E^{\epsilon}\times B^{\epsilon})=0,
\end{equation}
where $\mathcal{L}=(I+\epsilon^{2}\Delta^{2})^{-1}$. It is easily
to see that $\mathcal{L}$ satisfies the following properties:
\begin{equation}\label{4.4}
\left\{\!\!
\begin{array}{lc}
(1)\ \|\mathcal{L}f\|_{H^{k}}\leq \|f\|_{H^{k}},\ \forall\ k\in
\mathbb{R},\\
(2)\ (\mathcal{L}f,f)=\int_{\mathbb{R}^{d}}(\mathcal{L}f)\cdot
\bar{f}dx\geq 0,\\
(3)\ (\mathcal{L}f,g)=(f,\mathcal{L}g),\\
(4)\ \mathcal{L}\ commutes\ with\ Fourier\ multipler\ such\ as\
\Lambda^{s}, \nabla\ and\ so\ on.
\end{array}
\right.
\end{equation}
Due to these properties, the operator $\mathcal{L}$ can be easily
dealt with in the following estimates.

Since \eqref{4.1} holds, we have
\begin{equation}\label{4.5}
\|E_{0}^{\epsilon}\|_{H^{s+1}}+\|n_{0}^{\epsilon}\|_{H^{s}}+\|n_{1}^{\epsilon}\|_{H^{s-1}\cap
\dot{H}^{-1}}+\|B_{0}^{\epsilon}\|_{H^{s}\cap
\dot{H}^{-1}}+\|B_{1}^{\epsilon}\|_{H^{s-2}\cap \dot{H}^{-2}}\leq
c_{0},
\end{equation}
where the magnitude of $c_{0}$   depends only on
$\|E_{0}\|_{H^{s+1}}$, $ \|n_{0}\|_{H^{s}}$,
$\|n_{1}\|_{H^{s-1}\cap \dot{H}^{-1}}$, $\|B_{0}\|_{H^{s}\cap
\dot{H}^{-1}}$ and $\|B_{1}\|_{H^{s-2}\cap \dot{H}^{-2}}$.\\

$\bullet$ \textbf{Low order norm estimates.}\\
By the conserved quantities \eqref{3.2} and \eqref{3.3}, there
holds
\begin{align}\label{4.6}
&\|E^{\epsilon}\|_{H^{1}}^{2}+\|n^{\epsilon}\|_{L^{2}}^{2}+\|n^{\epsilon}_{t}\|_{\dot{H}^{-1}}^{2}
+\|B^{\epsilon}\|_{L^{2}}^{2}+\|B^{\epsilon}\|_{\dot{H}^{-1}}^{2}+\|B^{\epsilon}_{t}\|_{\dot{H}^{-2}}^{2}\nonumber\\
&\quad \leq
C(c_{0})+\left|\int_{\mathbb{R}^{d}}n^{\epsilon}|E^{\epsilon}|^{2}dx\right|
+\left|\int_{\mathbb{R}^{d}}(E^{\epsilon}\times
\overline{E^{\epsilon}})\cdot B^{\epsilon}dx\right|.
\end{align}
Then by Cauchy-Schwarz inequality, we have
\begin{align}\label{4.7}
\begin{split}
&\left|\int_{\mathbb{R}^{d}}n^{\epsilon}|E^{\epsilon}|^{2}dx\right|\leq
\|n^{\epsilon}\|_{L^{2}}\|E^{\epsilon}\|_{L^{4}}^{2}\leq
\frac{1}{2}\|n^{\epsilon}\|_{L^{2}}^{2}+\frac{1}{2}\|E^{\epsilon}\|_{L^{4}}^{4},\\
& \left|\int_{\mathbb{R}^{d}}(E^{\epsilon}\times
\overline{E^{\epsilon}})\cdot B^{\epsilon}dx\right|\leq
\frac{1}{2}\|B^{\epsilon}\|_{L^{2}}^{2}+\frac{1}{2}\|E^{\epsilon}\|_{L^{4}}^{4}.
\end{split}
\end{align}

Since $|E^{\epsilon}|^{2}_{t}=2\mathrm{Im}(iE^{\epsilon}_{t}\cdot
\overline{E^{\epsilon}})$, then \eqref{4.3} yields
\begin{equation}\label{4.8}
|E^{\epsilon}|^{2}_{t}=2\mathrm{Im}\Big[\big(-\mathcal{L}\nabla(\nabla\cdot
E^{\epsilon})+\alpha\mathcal{L}\nabla\times(\nabla\times
E^{\epsilon})+\mathcal{L}(n^{\epsilon}E^{\epsilon})-i\mathcal{L}(
E^{\epsilon}\times
B^{\epsilon})\big)\cdot\overline{E^{\epsilon}}\Big],
\end{equation}
hence we can obtain
\begin{align*}
\frac{d}{dt}\int_{\mathbb{R}^{d}}|E^{\epsilon}|^{4}dx&=2\int_{\mathbb{R}^{d}}|E^{\epsilon}|^{2}|E^{\epsilon}|^{2}_{t}dx\\
&\leq
C(\|E^{\epsilon}\|_{H^{s+1}}^{4}+\|E^{\epsilon}\|_{H^{s+1}}^{4}(\|n^{\epsilon}\|_{H^{s}}+\|B^{\epsilon}\|_{H^{s}}))\\
&\leq
C(\|E^{\epsilon}\|_{H^{s+1}}^{2}+\|n^{\epsilon}\|_{H^{s}}^{2}+\|B^{\epsilon}\|_{H^{s}}^{2}+1)^{3}
\end{align*}
which implies that
\begin{equation}\label{4.9}
\|E^{\epsilon}\|_{L^{4}}^{4}\leq
C+C\int_{0}^{t}(\|E^{\epsilon}\|_{H^{s+1}}^{2}+\|n^{\epsilon}\|_{H^{s}}^{2}+\|B^{\epsilon}\|_{H^{s}}^{2}+1)^{3}d\tau.
\end{equation}
Putting \eqref{4.6}, \eqref{4.7}, \eqref{4.9} together, then we
have
\begin{equation}\label{4.10}
\mathrm{LHS\ of}\ \eqref{4.6}\leq
C+C\int_{0}^{t}(\|E^{\epsilon}\|_{H^{s+1}}^{2}+\|n^{\epsilon}\|_{H^{s}}^{2}+\|B^{\epsilon}\|_{H^{s}}^{2}+1)^{3}d\tau.
\end{equation}

$\bullet$ \textbf{High order norm estimates.}\\

Applying the operator $\Lambda^{s}$ to equation \eqref{4.3}, then
one has
\begin{equation*}
i\Lambda^{s}E_{t}^{\epsilon}+\Lambda^{s}\mathcal{L}\nabla(\nabla\cdot
E^{\epsilon})-\alpha\Lambda^{s}\mathcal{L}\nabla\times(\nabla\times
E^{\epsilon})-\Lambda^{s}\mathcal{L}(n^{\epsilon}E^{\epsilon})+i\Lambda^{s}\mathcal{L}(
E^{\epsilon}\times B^{\epsilon})=0.
\end{equation*}
Taking inner product of this equation with
$-\Lambda^{s}\nabla(\nabla\cdot
E^{\epsilon})+\alpha\Lambda^{s}\nabla\times(\nabla\times
E^{\epsilon})$, and using the properties for $\mathcal{L}$ given
by \eqref{4.4}, one can obtain
\begin{align}
&\frac{d}{dt}(\|\Lambda^{s}(\nabla\cdot
E^{\epsilon})\|_{L^{2}}^{2}+\alpha\|\Lambda^{s}(\nabla\times
E^{\epsilon})\|_{L^{2}}^{2})\nonumber\\
&\quad
=2\mathrm{Im}\int_{\mathbb{R}^{d}}\mathcal{L}\Lambda^{s}(n^{\epsilon}E^{\epsilon})\cdot
(-\Lambda^{s}\nabla(\nabla\cdot
\overline{E^{\epsilon}})+\alpha\Lambda^{s}\nabla\times(\nabla\times
\overline{E^{\epsilon}}))dx\nonumber\\
&\quad\quad+ 2\mathrm{Im}\
i\int_{\mathbb{R}^{d}}\mathcal{L}\Lambda^{s}(E^{\epsilon}\times
B^{\epsilon})\cdot (\Lambda^{s}\nabla(\nabla\cdot
\overline{E^{\epsilon}})-\alpha\Lambda^{s}\nabla\times(\nabla\times
\overline{E^{\epsilon}}))dx\nonumber\\
&\quad =:I_{1}+I_{2}.\label{4.11}
\end{align}

We multiply equation \eqref{3.1b} by
$\Lambda^{2s-2}n^{\epsilon}_{t}$ and get
\begin{align}
&\frac{1}{2}\frac{d}{dt}\left(\|\Lambda^{s-1}n^{\epsilon}_{t}\|_{L^{2}}^{2}+\|\Lambda^{s}n^{\epsilon}\|_{L^{2}}^{2}
+2\int_{\mathbb{R}^{d}}\Lambda^{s}|E^{\epsilon}|^{2}\Lambda^{s}n^{\epsilon}dx\right)\nonumber\\
&\quad=\int_{\mathbb{R}^{d}}\Lambda^{s}n^{\epsilon}\Lambda^{s}|E^{\epsilon}|^{2}_{t}dx\nonumber\\
&\quad =2\mathrm{Im}\int_{\mathbb{R}^{d}}\Lambda^{s}n^{\epsilon}
\Lambda^{s}\Big[\big(-\mathcal{L}\nabla(\nabla\cdot
E^{\epsilon})+\alpha\mathcal{L}\nabla\times(\nabla\times
E^{\epsilon})\big)\cdot\overline{E^{\epsilon}}\Big]dx\nonumber\\
&\quad
\quad+2\mathrm{Im}\int_{\mathbb{R}^{d}}\Lambda^{s}n^{\epsilon}
\Lambda^{s}\Big[\big(\mathcal{L}(n^{\epsilon}E^{\epsilon})-i\mathcal{L}(
E^{\epsilon}\times B^{\epsilon})\big)\cdot\overline{E^{\epsilon}}\Big]dx\nonumber\\
&\quad =:I_{3}+I_{4}.\label{4.12}
\end{align}
It is obvious that
\begin{align}\label{4.13}
|I_{4}|&\leq
C\|n^{\epsilon}\|_{H^{s}}(\|n^{\epsilon}\|_{H^{s}}\|E^{\epsilon}\|_{H^{s+1}}^{2}+\|B^{\epsilon}\|_{H^{s}}\|E^{\epsilon}\|_{H^{s+1}}^{2})\nonumber\\
&\leq
C(\|n^{\epsilon}\|_{H^{s}}^{2}+\|B^{\epsilon}\|_{H^{s}}^{2}+\|E^{\epsilon}\|_{H^{s+1}}^{2})^{2}.
\end{align}
Now we estimate $I_{1}+I_{3}$. Since
\begin{align*}
I_{1}+I_{3}&=2\mathrm{Im}\int_{\mathbb{R}^{d}}\Lambda^{s+1}(n^{\epsilon}E^{\epsilon})\cdot
\big(-\Lambda^{s-1}\mathcal{L}\nabla(\nabla\cdot
\overline{E^{\epsilon}})+\Lambda^{s-1}\alpha\mathcal{L}\nabla\times(\nabla\times
\overline{E^{\epsilon}})\big)dx\\
&\quad +2\mathrm{Im}\int_{\mathbb{R}^{d}}\Lambda^{s}n^{\epsilon}
\Lambda^{s}\Big[\big(-\mathcal{L}\nabla(\nabla\cdot
E^{\epsilon})+\alpha\mathcal{L}\nabla\times(\nabla\times
E^{\epsilon})\big)\cdot\overline{E^{\epsilon}}\Big]dx,
\end{align*}
if we take $f=E^{\epsilon}$,
$g=-\Lambda^{-2}\mathcal{L}\nabla(\nabla\cdot
E^{\epsilon})+\Lambda^{-2}\alpha\mathcal{L}\nabla\times(\nabla\times
E^{\epsilon})$, $h=n^{\epsilon}$, then Lemma \ref{lem4.1} below
\begin{align}\label{4.14}
I_{1}+I_{3}&=2\mathrm{Im}\int_{\mathbb{R}^{d}}\Lambda^{s+1}(fh)\Lambda^{s+1}\bar{g}dx
+2\mathrm{Im}\int_{\mathbb{R}^{d}}\Lambda^{s}(\bar{f}\Lambda^{2}g)\Lambda^{s}hdx\nonumber\\
&\leq C\|f\|_{H^{s+1}}\|g\|_{H^{s+1}}\|h\|_{H^{s}}\ \nonumber\\
&\leq C
\|E^{\epsilon}\|_{H^{s+1}}^{2}\|n^{\epsilon}\|_{H^{s}}\nonumber\\
&\leq
C(\|E^{\epsilon}\|_{H^{s+1}}^{2}+\|n^{\epsilon}\|_{H^{s}}^{2}+1)^{2}.
\end{align}

Multiplying equation \eqref{3.1c} by
$\Lambda^{2s-4}B^{\epsilon}_{t}$, and using the fact
\begin{align}
(E^{\epsilon}\times
&\overline{E^{\epsilon}})_{t}=2i\mathrm{Im}(E^{\epsilon}_{t}\times
\overline{E^{\epsilon}})=-2i\mathrm{Re}(iE^{\epsilon}_{t}\times
\overline{E^{\epsilon}})\nonumber\\
& =2i \mathrm{Re}\Big[\big(\mathcal{L}\nabla(\nabla\cdot
E^{\epsilon})-\alpha\mathcal{L}\nabla\times(\nabla\times
E^{\epsilon})-\mathcal{L}(n^{\epsilon}E^{\epsilon})+i\mathcal{L}(
E^{\epsilon}\times B^{\epsilon})\big)\times
\overline{E^{\epsilon}}\Big],\label{4.15}
\end{align}
then we obtain
\begin{align}
&\frac{1}{2}\frac{d}{dt}\left(\|\Lambda^{s-2}B^{\epsilon}_{t}\|_{L^{2}}^{2}+\|\Lambda^{s}B^{\epsilon}\|_{L^{2}}^{2}
+\|\Lambda^{s-1}B^{\epsilon}\|_{L^{2}}^{2}
+2i\int_{\mathbb{R}^{d}}\Lambda^{s}(E^{\epsilon}\times \overline{E^{\epsilon}}) \Lambda^{s}B^{\epsilon}dx\right)\nonumber\\
&\quad=i\int_{\mathbb{R}^{d}}\Lambda^{s}B^{\epsilon}\Lambda^{s}(E^{\epsilon}\times \overline{E^{\epsilon}})_{t}dx\nonumber\\
&\quad =-2\mathrm{Re}\int_{\mathbb{R}^{d}}\Lambda^{s}B^{\epsilon}
\Lambda^{s}\Big[\big(\mathcal{L}\nabla(\nabla\cdot
E^{\epsilon})-\alpha\mathcal{L}\nabla\times(\nabla\times
E^{\epsilon})\big)\times\overline{E^{\epsilon}}\Big]dx\nonumber\\
&\quad
\quad-2\mathrm{Re}\int_{\mathbb{R}^{d}}\Lambda^{s}B^{\epsilon}
\Lambda^{s}\Big[\big(-\mathcal{L}(n^{\epsilon}E^{\epsilon})+i\mathcal{L}(
E^{\epsilon}\times B^{\epsilon})\big)\times\overline{E^{\epsilon}}\Big]dx\nonumber\\
&\quad =:I_{5}+I_{6}.\label{4.16}
\end{align}
Again the term $I_{6}$ can be estimated easily
\begin{align}\label{4.17}
|I_{6}|&\leq
C\|B^{\epsilon}\|_{H^{s}}(\|n^{\epsilon}\|_{H^{s}}\|E^{\epsilon}\|_{H^{s+1}}^{2}
+\|B^{\epsilon}\|_{H^{s}}\|E^{\epsilon}\|_{H^{s+1}}^{2})\nonumber\\
&\leq
C(\|n^{\epsilon}\|_{H^{s}}^{2}+\|B^{\epsilon}\|_{H^{s}}^{2}+\|E^{\epsilon}\|_{H^{s+1}}^{2})^{2}.
\end{align}
We need to estimate $I_{2}+I_{5}$. Rewrite $I_{2}+I_{5}$ in the
form
\begin{align*}
I_{2}+I_{5}&=-2\mathrm{Re}\int_{\mathbb{R}^{d}}\Lambda^{s+1}(E^{\epsilon}\times
B^{\epsilon})\cdot \big(\Lambda^{s-1}\mathcal{L}\nabla(\nabla\cdot
\overline{E^{\epsilon}})-\Lambda^{s-1}\alpha\mathcal{L}\nabla\times(\nabla\times
\overline{E^{\epsilon}})\big)dx\\
&\quad +2\mathrm{Re}\int_{\mathbb{R}^{d}}\Lambda^{s}B^{\epsilon}
\Lambda^{s}\Big[\overline{E^{\epsilon}}\times\big(\mathcal{L}\nabla(\nabla\cdot
E^{\epsilon})-\alpha\mathcal{L}\nabla\times(\nabla\times
E^{\epsilon})\big)\Big]dx,
\end{align*}
if we again take $f=E^{\epsilon}$,
$g=\Lambda^{-2}\mathcal{L}\nabla(\nabla\cdot
E^{\epsilon})-\Lambda^{-2}\alpha\mathcal{L}\nabla\times(\nabla\times
E^{\epsilon})$, $\tilde{h}=B^{\epsilon}$, then by Lemma
\ref{lem4.1} below, we have
\begin{align}\label{4.18}
I_{2}+I_{5}&=-2\mathrm{Re}\int_{\mathbb{R}^{d}}\Lambda^{s+1}(f\times
\tilde{h})\Lambda^{s+1}\bar{g}dx
+2\mathrm{Re}\int_{\mathbb{R}^{d}}\Lambda^{s}(\bar{f}\times\Lambda^{2}g)\Lambda^{s}\tilde{h}dx\nonumber\\
&\leq C\|f\|_{H^{s+1}}\|g\|_{H^{s+1}}\|\tilde{h}\|_{H^{s}}\nonumber\\
&\leq C
\|E^{\epsilon}\|_{H^{s+1}}^{2}\|B^{\epsilon}\|_{H^{s}}\nonumber\\
&\leq
C(\|E^{\epsilon}\|_{H^{s+1}}^{2}+\|B^{\epsilon}\|_{H^{s}}^{2}+1)^{2}.
\end{align}

Combining \eqref{4.11}-\eqref{4.14} and \eqref{4.16}-\eqref{4.18}
together, we arrive at
\begin{align}\label{4.19}
\begin{split}
& \|\Lambda^{s}(\nabla\cdot
  E^{\epsilon})\|_{L^{2}}^{2}+\alpha\|\Lambda^{s}(\nabla\times
  E^{\epsilon})\|_{L^{2}}^{2}+\|\Lambda^{s-1}n^{\epsilon}_{t}\|_{L^{2}}^{2}+\|\Lambda^{s}n^{\epsilon}\|_{L^{2}}^{2}\\
 &\quad \quad
 +\|\Lambda^{s-2}B^{\epsilon}_{t}\|_{L^{2}}^{2}+\|\Lambda^{s}B^{\epsilon}\|_{L^{2}}^{2}+\|\Lambda^{s-1}B^{\epsilon}\|_{L^{2}}^{2}\\
 &\quad\leq
 C+C\int_{0}^{t}(\|n^{\epsilon}\|_{H^{s}}^{2}+\|B^{\epsilon}\|_{H^{s}}^{2}+\|E^{\epsilon}\|_{H^{s+1}}^{2})^{2}d\tau\\
 &\quad\quad
 +2\left|\int_{\mathbb{R}^{d}}\Lambda^{s}|E^{\epsilon}|^{2}\Lambda^{s}n^{\epsilon}dx\right|
 +2\left|\int_{\mathbb{R}^{d}}\Lambda^{s}(E^{\epsilon}\times \overline{E^{\epsilon}})
 \Lambda^{s}B^{\epsilon}dx\right|.
\end{split}
\end{align}
Using Cauchy-Schwarz inequality, we obtain
\begin{align}\label{4.20}
\left|\int_{\mathbb{R}^{d}}\Lambda^{s}|E^{\epsilon}|^{2}\Lambda^{s}n^{\epsilon}dx\right|\leq
\|\Lambda^{s}n^{\epsilon}\|_{L^{2}}\|\Lambda^{s}|E^{\epsilon}|^{2}\|_{L^{2}}\leq
\frac{1}{2}\|\Lambda^{s}n^{\epsilon}\|_{L^{2}}^{2}+\frac{1}{2}\|\Lambda^{s}|E^{\epsilon}|^{2}\|_{L^{2}}^{2}
\end{align}
and
\begin{align}\label{4.21}
\left|\int_{\mathbb{R}^{d}}\Lambda^{s}(E^{\epsilon}\times
\overline{E^{\epsilon}})
 \Lambda^{s}B^{\epsilon}dx\right|& \leq \|\Lambda^{s}B^{\epsilon}\|_{L^{2}}\|\Lambda^{s}(E^{\epsilon}\times \overline{E^{\epsilon}})
 \|_{L^{2}}\nonumber\\
 & \leq \frac{1}{2}\|\Lambda^{s}B^{\epsilon}\|_{L^{2}}^{2}+\frac{1}{2}\|\Lambda^{s}(E^{\epsilon}\times
 \overline{E^{\epsilon}})\|_{L^{2}}^{2}.
\end{align}
Then by \eqref{4.8} we have
\begin{align*}
\frac{d}{dt}\int_{\mathbb{R}^{d}}(\Lambda^{s}|E^{\epsilon}|^{2})^{2}dx
&=2\int_{\mathbb{R}^{d}}\Lambda^{s}|E^{\epsilon}|^{2}\Lambda^{s}(|E^{\epsilon}|^{2}_{t})dx\\
&\leq
C\big(\|E^{\epsilon}\|_{H^{s+1}}^{4}+\|E^{\epsilon}\|_{H^{s+1}}^{4}(\|n^{\epsilon}\|_{H^{s}}+\|B^{\epsilon}\|_{H^{s}})\big)\\
&\leq
C(\|E^{\epsilon}\|_{H^{s+1}}^{2}+\|n^{\epsilon}\|_{H^{s}}^{2}+\|B^{\epsilon}\|_{H^{s}}^{2}+1)^{3}
\end{align*}
which gives that
\begin{equation}\label{4.22}
\|\Lambda^{s}|E^{\epsilon}|^{2}\|_{L^{2}}^{2}\leq
C+C\int_{0}^{t}(\|E^{\epsilon}\|_{H^{s+1}}^{2}+\|n^{\epsilon}\|_{H^{s}}^{2}+\|B^{\epsilon}\|_{H^{s}}^{2}+1)^{3}d\tau.
\end{equation}
Using \eqref{4.15} and the same argument, we can also get
\begin{equation}\label{4.23}
\|\Lambda^{s}(E^{\epsilon}\times
\overline{E^{\epsilon}})\|_{L^{2}}^{2}\leq
C+C\int_{0}^{t}(\|E^{\epsilon}\|_{H^{s+1}}^{2}+\|n^{\epsilon}\|_{H^{s}}^{2}+\|B^{\epsilon}\|_{H^{s}}^{2}+1)^{3}d\tau.
\end{equation}
Inserting \eqref{4.20}-\eqref{4.23} into \eqref{4.19}, hence we
have
\begin{equation}\label{4.24}
\mathrm{LHS\ of}\ \eqref{4.19}\leq
C+C\int_{0}^{t}(\|E^{\epsilon}\|_{H^{s+1}}^{2}+\|n^{\epsilon}\|_{H^{s}}^{2}+\|B^{\epsilon}\|_{H^{s}}^{2}+1)^{3}d\tau.
\end{equation}

$\bullet$ \textbf{Conclusions.}\\

It concludes from the low order estimate \eqref{4.10} and the high
order estimate \eqref{4.24} that
\begin{align*}
&\|E^{\epsilon}\|_{H^{s+1}}^{2}+\|n^{\epsilon}\|_{H^{s}}^{2}+\|B^{\epsilon}\|_{H^{s}\cap
\dot{H}^{-1}}^{2} +\|n^{\epsilon}_{t}\|_{H^{s-1}\cap
\dot{H}^{-1}}^{2}+\|B^{\epsilon}_{t}\|_{H^{s-2}\cap
\dot{H}^{-2}}^{2}\\
&\quad\quad\leq
C+C\int_{0}^{t}(\|E^{\epsilon}\|_{H^{s+1}}^{2}+\|n^{\epsilon}\|_{H^{s}}^{2}+\|B^{\epsilon}\|_{H^{s}}^{2}+1)^{3}d\tau,
\end{align*}
where $C$ depends on $c_{0}$, hence by Lemma \ref{lem4.2} below,
we know there exist $T>0$ and $C>0$ both independent of $\epsilon$
such that for all $t\in [0,T]$
\begin{equation}\label{4.25}
\|E^{\epsilon}\|_{H^{s+1}}^{2}+\|n^{\epsilon}\|_{H^{s}}^{2}+\|B^{\epsilon}\|_{H^{s}\cap
\dot{H}^{-1}}^{2} +\|n^{\epsilon}_{t}\|_{H^{s-1}\cap
\dot{H}^{-1}}^{2}+\|B^{\epsilon}_{t}\|_{H^{s-2}\cap
\dot{H}^{-2}}^{2}\leq C,
\end{equation}
from which Proposition \ref{prop4.1} follows.
\end{proof}

Now we are going to prove the following lemma which is used in the
proof of Proposition \ref{prop4.1}.
\begin{lemma}\label{lem4.1}
Assume that $f,g\in H^{s+1}(\mathbb{R}^{d})$ are $\mathbb{C}^{3}$
valued functions, and $\tilde{h}\in H^{s}(\mathbb{R}^{d})$ is a
$\mathbb{R}^{3}$ valued function, and $h\in H^{s}(\mathbb{R}^{d})$
is a real valued function, $s>\frac{d}{2}$. Then the following
three estimates hold:
\begin{align}
&\left|\mathrm{Im}\int_{\mathbb{R}^{d}}\Lambda^{s+1}(fh)\cdot\Lambda^{s+1}\bar{g}dx
+\mathrm{Im}\int_{\mathbb{R}^{d}}\Lambda^{s}(\bar{f}\cdot\Lambda^{2}g)\Lambda^{s}hdx\right|\leq
C\|f\|_{H^{s+1}}\|g\|_{H^{s+1}}\|h\|_{H^{s}},\label{4.26}\\
&\left|\mathrm{Re}\int_{\mathbb{R}^{d}}\Lambda^{s+1}(fh)\cdot\Lambda^{s+1}\bar{g}dx
-\mathrm{Re}\int_{\mathbb{R}^{d}}\Lambda^{s}(\bar{f}\cdot\Lambda^{2}g)\Lambda^{s}hdx\right|\leq
C\|f\|_{H^{s+1}}\|g\|_{H^{s+1}}\|h\|_{H^{s}},\label{4.27}\\
&\left|\mathrm{Re}\int_{\mathbb{R}^{d}}\Lambda^{s+1}(f\times
\tilde{h})\cdot\Lambda^{s+1}\bar{g}dx
-\mathrm{Re}\int_{\mathbb{R}^{d}}\Lambda^{s}(\bar{f}\times\Lambda^{2}g)\cdot\Lambda^{s}\tilde{h}dx\right|\leq
C\|f\|_{H^{s+1}}\|g\|_{H^{s+1}}\|\tilde{h}\|_{H^{s}}.\label{4.28}
\end{align}
\end{lemma}

\begin{proof} We first show \eqref{4.26}. Denote the LHS of
\eqref{4.26} by $|J|=|J_{1}+J_{2}|$. The term $J_{1}$ can be
written as
\begin{align}
J_{1}&=\mathrm{Im}\int_{\mathbb{R}^{d}}[\Lambda^{s+1}(fh)-f\Lambda^{s+1}h]\cdot\Lambda^{s+1}\bar{g}dx
+\mathrm{Im}\int_{\mathbb{R}^{d}}f\Lambda^{s+1}h\cdot\Lambda^{s+1}\bar{g}dx\nonumber\\
&=\mathrm{Im}\int_{\mathbb{R}^{d}}[\Lambda^{s+1}(fh)-f\Lambda^{s+1}h]\cdot\Lambda^{s+1}\bar{g}dx\nonumber\\
&\quad+\mathrm{Im}\int_{\mathbb{R}^{d}}\Lambda^{s}h[\Lambda(f\cdot\Lambda^{s+1}\bar{g})-f\cdot\Lambda^{s+2}\bar{g}]dx
+\mathrm{Im}\int_{\mathbb{R}^{d}}f\Lambda^{s}h\cdot\Lambda^{s+2}\bar{g}dx\nonumber\\
&=:J_{11}+J_{12}+J_{13}.\label{4.29}
\end{align}
From the commutator estimate \eqref{3.8}, we have
\begin{align}
|J_{11}|&\leq C(\|\nabla
f\|_{L^{\infty}}\|h\|_{H^{s}}+\|f\|_{H^{s+1}}\|h\|_{L^{\infty}})\|g\|_{H^{s+1}}\nonumber\\
&\leq C\|f\|_{H^{s+1}}\|g\|_{H^{s+1}}\|h\|_{H^{s}},\label{4.30}
\end{align}
since $s>\frac{d}{2}$. The term $J_{12}$ can be estimated by
commutator estimate \eqref{3.8} again
\begin{equation*}
|J_{12}|\leq C\|h\|_{H^{s}}(\|\nabla
f\|_{L^{\infty}}\|g\|_{H^{s+1}}+\|f\|_{H^{1,p}}\|g\|_{H^{s+1,q}}),
\end{equation*}
where $\frac{1}{p}+\frac{1}{q}=\frac{1}{2}$.  Since
$s>\frac{d}{2}$, we have $H^{s+1}\hookrightarrow H^{1,p}$ for all
$p\in [2,\infty]$, then
\begin{equation*}
|J_{12}|\leq
C\|h\|_{H^{s}}(\|f\|_{H^{s+1}}\|g\|_{H^{s+1}}+\|f\|_{H^{s+1}}\|g\|_{H^{s+1,q}}),
\end{equation*}
letting $q\rightarrow 2^{+}$ in the above inequality, then we get
\begin{equation}\label{4.31}
|J_{12}|\leq C\|f\|_{H^{s+1}}\|g\|_{H^{s+1}}\|h\|_{H^{s}}.
\end{equation}

Now we estimate $J_{2}$. It is obvious that
\begin{align}
J_{2}&=\mathrm{Im}\int_{\mathbb{R}^{d}}[\Lambda^{s}(\bar{f}\cdot\Lambda^{2}g)-\bar{f}\cdot\Lambda^{s+2}g]\Lambda^{s}hdx
+\mathrm{Im}\int_{\mathbb{R}^{d}}\bar{f}\cdot\Lambda^{s+2}g\Lambda^{s}hdx\nonumber\\
&=:J_{21}+J_{22}.\label{4.31}
\end{align}
Using commutator estimate \eqref{3.8}, we obtain
\begin{equation*}
|J_{21}|\leq C(\|\nabla
f\|_{L^{\infty}}\|g\|_{H^{s+1}}+\|f\|_{H^{s,p}}\|g\|_{H^{2,q}})\|h\|_{H^{s}},
\end{equation*}
where we select $p$, $q$ satisfying
$$
\frac{1}{p}+\frac{1}{q}=\frac{1}{2},\ \frac{d}{2}\leq
1+\frac{d}{p},\ 1+\frac{d}{2}-\frac{d}{q}\leq s,
$$
hence we have $H^{s+1}\hookrightarrow H^{s,p}$ and
$H^{s+1}\hookrightarrow H^{2,q}$,  and we therefore get
\begin{equation}\label{4.33}
|J_{21}|\leq C\| f\|_{H^{s+1}}\|g\|_{H^{s+1}}\|h\|_{H^{s}}.
\end{equation}

Since $J_{13}+J_{22}=0$, now \eqref{4.26} follows from
\eqref{4.29}-\eqref{4.33}. If we replace $f$ by $i\cdot f$, then
\eqref{4.27} is a direct consequence of \eqref{4.26}. Moreover, if
we expand the term $\Lambda^{s+1}(f\times
\tilde{h})\cdot\Lambda^{s+1}\bar{g}$ and
$\Lambda^{s}(\bar{f}\times\Lambda^{2}g)\cdot\Lambda^{s}\tilde{h}$
by the definition of dot product and cross product, then one can
easily see that \eqref{4.28} reduces to \eqref{4.27}. We thus
finish the proof of Lemma \ref{lem4.1}.
\end{proof}

Now we end this section with the following elementary lemma.

\begin{lemma}\label{lem4.2}
Let $u(t)$ be a continuous and nonnegative function defined on
$\mathbb{R}^{+}$, and suppose $u$ obeys the integral inequality
\begin{equation}\label{4.34}
u(t)\leq A_{1}+A_{2}\int_{0}^{t}u^{\kappa}(\tau)d\tau,\ \kappa> 1
\end{equation}
for all $t\geq 0$, where $A_{1}$, $A_{2}> 0$. Moreover, let
\begin{align*}
v(t)=\frac{A_{1}}{(1-(\kappa-1)A_{1}^{\kappa-1}A_{2}t)^{1/(\kappa-1)}},
\end{align*}
and
\begin{align*}
T^{*}=\frac{1}{(\kappa-1)A_{1}^{\kappa-1}A_{2}}.
\end{align*} Then
$$
u(t)\leq v(t),\ \forall\ t\in[0,T^{*}).
$$
\end{lemma}
\begin{proof}It is obviously that $v$ satisfies
\begin{equation*}
v(t)= A_{1}+A_{2}\int_{0}^{t}v^{\kappa}(\tau)d\tau,
\end{equation*}
hence this equality and \eqref{4.34} yield
$$
u(t)-v(t)\leq A_{2}\int_{0}^{t}\kappa(\theta
u(\tau)+(1-\theta)v(\tau))^{\kappa-1}(u(\tau)-v(\tau))d\tau,
$$
which implies that
$$
w(t)\leq A_{2}\int_{0}^{t}\kappa(\theta
u(\tau)+(1-\theta)v(\tau))^{\kappa-1}w(\tau)d\tau,
$$
where $w(t):=(u(t)-v(t))_{+}=\max\{u(t)-v(t),0\}$. Fix any
$T<T^{*}$, since both $u$ and $v$ are nonnegative and continuous,
we have $0\leq u(t)$, $v(t)\leq K$ for all $t\in [0,T]$. Thus we
obtain
$$
w(t)\leq \kappa A_{2}K^{\kappa-1}\int_{0}^{t}w(\tau)d\tau,
$$
by Gronwall's inequality, then we have $w(t)=0$ for all $t\in
[0,T]$. Let $T\rightarrow T^{*}$, Lemma \ref{lem4.2} thus follows.
\end{proof}
%
%
%
%
%

\section{Strong convergence of the approximate solutions}
\setcounter{section}{5}\setcounter{equation}{0}

Under the prior estimates given in Proposition \ref{prop4.1}, we
now show that the solutions
$(E^{\epsilon},n^{\epsilon},B^{\epsilon})$ to the regularized
system \eqref{3.1a}-\eqref{3.1d} form a Cauchy sequence in the low
order norm $C\big([0,T);H^{1}\oplus L^{2}\oplus (L^{2}\cap
\dot{H}^{-1})\big)$. Namely, we are going to prove the following
lemma.

\begin{lemma}\label{lem5.1}
With the same assumptions as Proposition \ref{prop4.1}, then the
family $(E^{\epsilon},n^{\epsilon},B^{\epsilon})$ forms a Cauchy
sequence in $C\big([0,T];H^{1}\oplus L^{2}\oplus (L^{2}\cap
\dot{H}^{-1})\big)$, a.e. there holds
\begin{equation}\label{5.1}
\sup\limits_{t\in[0,T]}(\|E^{\epsilon}-E^{\epsilon'}\|_{H^{1}}+\|n^{\epsilon}-n^{\epsilon'}\|_{L^{2}}
+\|B^{\epsilon}-B^{\epsilon'}\|_{L^{2}\cap
\dot{H}^{-1}})\rightarrow 0,\ \epsilon,\ \epsilon'\rightarrow 0.
\end{equation}
Moreover, there holds
\begin{equation}\label{A5.1}
\sup\limits_{t\in[0,T]}(\|E^{\epsilon}_{t}-E^{\epsilon'}_{t}\|_{H^{-1}}+\|n^{\epsilon}_{t}-n^{\epsilon'}_{t}\|_{\dot{H}^{-1}}
+\|B^{\epsilon}_{t}-B^{\epsilon'}_{t}\|_{\dot{H}^{-2}})\rightarrow
0,\ \epsilon,\ \epsilon'\rightarrow 0.
\end{equation}
\end{lemma}

\begin{proof} For brevity, we set $E=E^{\epsilon}-E^{\epsilon'}$,
$n=n^{\epsilon}-n^{\epsilon'}$, $B=B^{\epsilon}-B^{\epsilon'}$.
Since $(E^{\epsilon},n^{\epsilon},B^{\epsilon})$ and
$(E^{\epsilon'},n^{\epsilon'},B^{\epsilon'})$ both satisfy the
regularized system \eqref{3.1a}-\eqref{3.1c}, then $(E,n,B)$
satisfies the equation
\begin{subequations}
\begin{align}
& iE_{t}+\mathcal{L}\nabla(\nabla\cdot
E)-\alpha\mathcal{L}\nabla\times(\nabla\times
E)-\mathcal{L}(nE^{\epsilon}+n^{\epsilon'}E)+i \mathcal{L}(E\times B^{\epsilon}+E^{\epsilon'}\times B)=0, \label{5.2a}\\
&n_{tt}-\Delta n=\Delta (|E^{\epsilon}|^{2}-|E^{\epsilon'}|^{2}),\label{5.2b}\\
&B_{tt}+\Delta ^{2}B-\Delta B=-i\Delta^{2}(E^{\epsilon}\times
\overline{E^{\epsilon}}-E^{\epsilon'}\times
\overline{E^{\epsilon'}})\label{5.2c}
 \end{align}
with initial data
\begin{align}\label{5.2d}
\begin{split}
&E(0)=E^{\epsilon}_{0}-E^{\epsilon'}_{0},\
n(0)=n^{\epsilon}_{0}-n^{\epsilon'}_{0},\
n_{t}(0)=n^{\epsilon}_{1}-n^{\epsilon'}_{1},\\
&B(0)=B^{\epsilon}_{0}-B^{\epsilon'}_{0},\
B_{t}(0)=B^{\epsilon}_{1}-B^{\epsilon'}_{1}.
\end{split}
\end{align}
\end{subequations}

From equation \eqref{5.2a}, one can obtain
\begin{align}\label{5.3}
\begin{split}
\frac{d}{dt}\|E\|_{L^{2}}^{2}&=2\mathrm{Im}\big(\mathcal{L}(nE^{\epsilon}+n^{\epsilon'}E)-i
\mathcal{L}(E\times B^{\epsilon}+E^{\epsilon'}\times B),E\big)\\
&\leq C(\|n\|_{L^{2}}\|E\|_{L^{2}}+\|E\|_{L^{2}}^{2}+
\|E\|_{L^{2}}\|B\|_{L^{2}})\\
&\leq C(\|n\|_{L^{2}}^{2}+\|E\|_{L^{2}}^{2}+\|B\|_{L^{2}}^{2}),
\end{split}
\end{align}
where we have used the prior estimate \eqref{4.2} in the first
inequality above.

If we multiply equation \eqref{5.2a} by $-\nabla(\nabla\cdot
\bar{E})+\alpha\nabla\times(\nabla\times \bar{E})$, and integrate
the imaginary part of the result, then we have
\begin{align}\label{5.4}
\begin{split}
&\frac{d}{dt}(\|\nabla\cdot E\|_{L^{2}}^{2}+\alpha\|\nabla\times
E\|_{L^{2}}^{2})\\
&\quad =2 \mathrm{Im}
\int_{\mathbb{R}^{d}}\mathcal{L}(nE^{\epsilon})\cdot
(-\nabla(\nabla\cdot \bar{E})+\alpha\nabla\times(\nabla\times
\bar{E}))dx\\
&\quad\quad +2 \mathrm{Im}
\int_{\mathbb{R}^{d}}\mathcal{L}(n^{\epsilon'}E)\cdot
(-\nabla(\nabla\cdot \bar{E})+\alpha\nabla\times(\nabla\times
\bar{E}))dx\\
&\quad\quad -2 \mathrm{Im}\ i
\int_{\mathbb{R}^{d}}\mathcal{L}(E\times B^{\epsilon})\cdot
(-\nabla(\nabla\cdot \bar{E})+\alpha\nabla\times(\nabla\times
\bar{E}))dx\\
&\quad\quad -2 \mathrm{Im}\ i
\int_{\mathbb{R}^{d}}\mathcal{L}(E^{\epsilon'}\times B)\cdot
(-\nabla(\nabla\cdot \bar{E})+\alpha\nabla\times(\nabla\times
\bar{E}))dx\\
& \quad =:K_{1}+K_{2}+K_{3}+K_{4}.
\end{split}
\end{align}
Using integrating by parts, Sobolev inequality, and the prior
estimate \eqref{4.2}, we can easily get
\begin{equation}\label{5.5}
|K_{2}|+|K_{3}|\leq C\|E\|_{H^{1}}^{2}.
\end{equation}

Now taking inner product to equation \eqref{5.2b} with
$\Lambda^{-2}n_{t}$, then
\begin{equation}\label{5.6}
\frac{1}{2}\frac{d}{dt}\left(\|\Lambda^{-1}n_{t}\|_{L^{2}}^{2}+\|n\|_{L^{2}}^{2}+2\int_{\mathbb{R}^{d}}(|E^{\epsilon}|^{2}-|E^{\epsilon'}|^{2})\cdot
n
dx\right)=\int_{\mathbb{R}^{d}}(|E^{\epsilon}|^{2}_{t}-|E^{\epsilon'}|^{2}_{t})ndx.
\end{equation}
From \eqref{4.8}, we compute
\begin{align}\label{5.7}
\begin{split}
|E^{\epsilon}|^{2}_{t}-|E^{\epsilon'}|^{2}_{t}&=2\mathrm{Im}\Big[\big(-\mathcal{L}\nabla(\nabla\cdot
E)  +\alpha\mathcal{L}\nabla\times(\nabla\times
E)\\
&\quad +\mathcal{L}(nE^{\epsilon}+n^{\epsilon'}E)-i\mathcal{L}(
E\times B^{\epsilon}+E^{\epsilon'}\times B) \big)\cdot
\overline{E^{\epsilon}}\Big]\\
&\quad +2\mathrm{Im}\Big[\big(-\mathcal{L}\nabla(\nabla\cdot
E^{\epsilon'})+\alpha\mathcal{L}\nabla\times(\nabla\times
E^{\epsilon'})\\
&\quad +\mathcal{L}(n^{\epsilon'}E^{\epsilon'})-i\mathcal{L}(
E^{\epsilon'}\times B^{\epsilon'})\big)\cdot\bar{E}\Big].
\end{split}
\end{align}
Plugging this equality into \eqref{5.6}, and again using the prior
estimate \eqref{4.2}, sometimes integrating by parts, and using
Sobolev inequality, then we have
\begin{align}\label{5.8}
\begin{split}
\mathrm{LHS\ of} \eqref{5.6} &\leq
2\mathrm{Im}\int_{\mathbb{R}^{d}}\Big[\big(-\mathcal{L}\nabla(\nabla\cdot
E)  +\alpha\mathcal{L}\nabla\times(\nabla\times E)\big)\cdot
\overline{E^{\epsilon}}\Big]ndx\\
&\quad +C (\|E\|_{H^{1}}^{2}+\|n\|_{L^{2}}^{2}+\|B\|_{L^{2}}^{2})\\
&\quad =: K_{5}+C
(\|E\|_{H^{1}}^{2}+\|n\|_{L^{2}}^{2}+\|B\|_{L^{2}}^{2}).
\end{split}
\end{align}

For the equation \eqref{5.2c}, we multiply it by
$\Lambda^{-4}B_{t}$ and obtain
\begin{align}\label{5.9}
\begin{split}
&\frac{1}{2}\frac{d}{dt}\Big(\|\Lambda^{-2}B_{t}\|_{L^{2}}^{2}+\|B\|_{L^{2}}^{2}+\|\Lambda^{-1}B\|_{L^{2}}^{2}\\
&\quad +2i\int_{\mathbb{R}^{d}}(E^{\epsilon}\times
\overline{E^{\epsilon}} -E^{\epsilon'}\times
\overline{E^{\epsilon'}})\cdot B
dx\Big)=i\int_{\mathbb{R}^{d}}(E^{\epsilon}\times
\overline{E^{\epsilon}} -E^{\epsilon'}\times
\overline{E^{\epsilon'}})_{t}\cdot Bdx.
\end{split}
\end{align}
Using \eqref{4.15}, we have
\begin{align}\label{5.10}
\begin{split}
(E^{\epsilon}\times \overline{E^{\epsilon}} -E^{\epsilon'}\times
\overline{E^{\epsilon'}})_{t}&=2i\mathrm{Re}\Big[\big(\mathcal{L}\nabla(\nabla\cdot
E)  -\alpha\mathcal{L}\nabla\times(\nabla\times
E)\\
&\quad -\mathcal{L}(nE^{\epsilon}+n^{\epsilon'}E)+i\mathcal{L}(
E\times B^{\epsilon}+E^{\epsilon'}\times B) \big)\times
\overline{E^{\epsilon'}}\Big]\\
&\quad +2i\mathrm{Re}\Big[\big(\mathcal{L}\nabla(\nabla\cdot
E^{\epsilon})-\alpha\mathcal{L}\nabla\times(\nabla\times
E^{\epsilon})\\
&\quad -\mathcal{L}(n^{\epsilon}E^{\epsilon})-i\mathcal{L}(
E^{\epsilon}\times B^{\epsilon})\big)\times\bar{E}\Big].
\end{split}
\end{align}
Then by the same reasonings that lead to \eqref{5.7}, we have
\begin{align}\label{5.11}
\begin{split}
\mathrm{LHS\ of} \eqref{5.8} &\leq
-2\mathrm{Re}\int_{\mathbb{R}^{d}}\Big[\big(\mathcal{L}\nabla(\nabla\cdot
E)  -\alpha\mathcal{L}\nabla\times(\nabla\times E)\big)\times
\overline{E^{\epsilon'}}\Big]\cdot Bdx\\
&\quad +C (\|E\|_{H^{1}}^{2}+\|n\|_{L^{2}}^{2}+\|B\|_{L^{2}}^{2})\\
&\quad =: K_{6}+C
(\|E\|_{H^{1}}^{2}+\|n\|_{L^{2}}^{2}+\|B\|_{L^{2}}^{2}).
\end{split}
\end{align}

Note that $K_{1}+K_{5}=0$, $K_{4}+K_{6}=0$, hence inequalities
\eqref{5.3}-\eqref{5.5}, \eqref{5.9} and \eqref{5.11} yield
\begin{align*}
&\frac{d}{dt}\Big(\|E\|_{L^{2}}^{2}+\|\nabla\cdot
E\|_{L^{2}}^{2}+\alpha\|\nabla\times
E\|_{L^{2}}^{2}+\frac{1}{2}\|\Lambda^{-1}n_{t}\|_{L^{2}}^{2}+\frac{1}{2}\|n\|_{L^{2}}^{2}\\
&\quad\quad+\frac{1}{2} \|\Lambda^{-2}B_{t}\|_{L^{2}}^{2}
+\frac{1}{2}\|B\|_{L^{2}}^{2}+\frac{1}{2}\|\Lambda^{-1}B\|_{L^{2}}^{2}\\
&\quad\quad+\int_{\mathbb{R}^{d}}(|E^{\epsilon}|^{2}-|E^{\epsilon'}|^{2})\cdot
n dx
 +i\int_{\mathbb{R}^{d}}(E^{\epsilon}\times
\overline{E^{\epsilon}} -E^{\epsilon'}\times
\overline{E^{\epsilon'}})\cdot B dx\Big)\\
&\quad\leq
C(\|E\|_{H^{1}}^{2}+\|n\|_{L^{2}}^{2}+\|B\|_{L^{2}}^{2}).
\end{align*}
By integrating this inequality, we obtain
\begin{align}\label{5.12}
\begin{split}
&\|E\|_{H^{1}}^{2}+\|n\|_{L^{2}}^{2}+\|B\|_{L^{2}\cap
\dot{H}^{-1}}^{2}+\|n_{t}\|_{\dot{H}^{-1}}^{2}+\|B_{t}\|_{\dot{H}^{-2}}^{2}\\
&\quad \leq
C(\|E(0)\|_{H^{1}}^{2}+\|n(0)\|_{L^{2}}^{2}+\|B(0)\|_{L^{2}\cap
\dot{H}^{-1}}^{2}+\|n_{t}(0)\|_{\dot{H}^{-1}}^{2}+\|B_{t}(0)\|_{\dot{H}^{-2}}^{2})\\
&\quad\quad+
C\left|\int_{\mathbb{R}^{d}}(|E^{\epsilon}|^{2}-|E^{\epsilon'}|^{2})\cdot
n dx\right|
 +C\left|\int_{\mathbb{R}^{d}}(E^{\epsilon}\times
\overline{E^{\epsilon}} -E^{\epsilon'}\times
\overline{E^{\epsilon'}})\cdot B dx\right|\\
&\quad\quad
+C\int_{0}^{t}(\|E\|_{H^{1}}^{2}+\|n\|_{L^{2}}^{2}+\|B\|_{L^{2}}^{2})d\tau.
\end{split}
\end{align}
Applying Cauchy-Schwarz inequality, we have
\begin{align}\label{5.13}
\begin{split}
&\left|\int_{\mathbb{R}^{d}}(|E^{\epsilon}|^{2}-|E^{\epsilon'}|^{2})\cdot
n dx\right|\leq
\frac{1}{2}\|n\|_{L^{2}}^{2}+\frac{1}{2}\||E^{\epsilon}|^{2}-|E^{\epsilon'}|^{2}\|_{L^{2}}^{2},\\
&\left|\int_{\mathbb{R}^{d}}(E^{\epsilon}\times
\overline{E^{\epsilon}} -E^{\epsilon'}\times
\overline{E^{\epsilon'}})\cdot B dx\right|\leq
\frac{1}{2}\|B\|_{L^{2}}^{2}+\frac{1}{2}\|E^{\epsilon}\times
\overline{E^{\epsilon}} -E^{\epsilon'}\times
\overline{E^{\epsilon'}}\|_{L^{2}}^{2}.
\end{split}
\end{align}
We deduce from \eqref{5.7} that
\begin{align*}
\frac{d}{dt}\||E^{\epsilon}|^{2}-|E^{\epsilon'}|^{2}\|_{L^{2}}^{2}
&=2\int_{\mathbb{R}^{d}}(|E^{\epsilon}|^{2}-|E^{\epsilon'}|^{2})\cdot(|E^{\epsilon}|^{2}_{t}-|E^{\epsilon'}|^{2}_{t})dx\\
&\leq C(\|E\|_{H^{1}}^{2}+\|n\|_{L^{2}}^{2}+\|B\|_{L^{2}}^{2}),
\end{align*}
thus one can obtain from this inequality
\begin{equation}\label{5.14}
\||E^{\epsilon}|^{2}-|E^{\epsilon'}|^{2}\|_{L^{2}}^{2} \leq
C\|E(0)\|_{H^{1}}^{2}+
C\int_{0}^{t}(\|E\|_{H^{1}}^{2}+\|n\|_{L^{2}}^{2}+\|B\|_{L^{2}}^{2})d\tau.
\end{equation}
A similar argument yields(we shall use \eqref{5.10} instead)
\begin{equation}\label{5.15}
\|E^{\epsilon}\times \overline{E^{\epsilon}} -E^{\epsilon'}\times
\overline{E^{\epsilon'}}\|_{L^{2}}^{2} \leq C\|E(0)\|_{H^{1}}^{2}+
C\int_{0}^{t}(\|E\|_{H^{1}}^{2}+\|n\|_{L^{2}}^{2}+\|B\|_{L^{2}}^{2})d\tau.
\end{equation}

Putting \eqref{5.12}-\eqref{5.15} together gives
\begin{align}
&\|E\|_{H^{1}}^{2}+\|n\|_{L^{2}}^{2}+\|B\|_{L^{2}\cap
\dot{H}^{-1}}^{2}+\|n_{t}\|_{\dot{H}^{-1}}^{2}+\|B_{t}\|_{\dot{H}^{-2}}^{2}\nonumber\\
&\quad \leq
C(\|E(0)\|_{H^{1}}^{2}+\|n(0)\|_{L^{2}}^{2}+\|B(0)\|_{L^{2}\cap
\dot{H}^{-1}}^{2}+\|n_{t}(0)\|_{\dot{H}^{-1}}^{2}+\|B_{t}(0)\|_{\dot{H}^{-2}}^{2})\nonumber\\
&\quad\quad
+C\int_{0}^{t}(\|E\|_{H^{1}}^{2}+\|n\|_{L^{2}}^{2}+\|B\|_{L^{2}}^{2})d\tau.\label{5.16}
\end{align}
Since \eqref{4.1} holds, by Gronwall's inequality, we thus deduce
from \eqref{5.16} that $(E^{\epsilon},n^{\epsilon},B^{\epsilon})$
is a Cauchy sequence in $C\big([0,T];H^{1}\oplus L^{2}\oplus
(L^{2}\cap \dot{H}^{-1})\big)$. Moreover, from \eqref{5.16} and
the equation \eqref{4.3}, we get \eqref{A5.1}. Then Lemma
\ref{lem5.1} follows.
\end{proof}

\begin{remark} From \eqref{4.2}, \eqref{A4.1},
\eqref{5.1}, \eqref{A5.1} and the following interpolation in
Sobolev spaces
\begin{align*}
\|f\|_{\dot{H}^{s_{0}}}\leq C
\|f\|_{\dot{H}^{s_{1}}}^{1-\theta}\|f\|_{\dot{H}^{s_{2}}}^{\theta},\
s_{1}<s_{2},\ s_{1}\leq s_{0}\leq s_{2},\
\theta=\frac{s_{0}-s_{1}}{s_{2}-s_{1}},
\end{align*}
we have as $\epsilon,\ \epsilon'\rightarrow 0$
\begin{align}
\begin{split}\label{5.17}
&\sup\limits_{t\in[0,T]}(\|E^{\epsilon}-E^{\epsilon'}\|_{H^{\tilde{s}+1}}+\|n^{\epsilon}-n^{\epsilon'}\|_{H^{\tilde{s}}}
+\|B^{\epsilon}-B^{\epsilon'}\|_{H^{\tilde{s}}\cap
\dot{H}^{-1}})\rightarrow 0,\\
&\sup\limits_{t\in[0,T]}(\|E^{\epsilon}_{t}-E^{\epsilon'}_{t}\|_{H^{\tilde{s}-1}}+\|n^{\epsilon}_{t}-n^{\epsilon'}_{t}\|_{H^{\tilde{s}-1}\cap\dot{H}^{-1}}
+\|B^{\epsilon}_{t}-B^{\epsilon'}_{t}\|_{H^{\tilde{s}-2}\cap\dot{H}^{-2}})\rightarrow
0,
\end{split}
\end{align}
for all $\tilde{s}<s$.
\end{remark}

\section{Proof of the main theorem}
\setcounter{section}{6}\setcounter{equation}{0}

\emph{Proof of Theorem \ref{thm1.1}.} For given $E_{0}\in
H^{s+1}$, $(n_{0},n_{1})\in H^{s}\oplus (H^{s-1}\cap
\dot{H}^{-1})$, $(B_{0},B_{1})\in (H^{s}\cap \dot{H}^{-1})\oplus
(H^{s-2}\cap \dot{H}^{-2})$, $s>\frac{d}{2}$, we choose
$(E_{0}^{\epsilon},n_{0}^{\epsilon},n_{1}^{\epsilon},B_{0}^{\epsilon},B_{1}^{\epsilon})$
sufficiently regular such that \eqref{4.1} holds. Then by the
strong convergence results \eqref{5.17}, we know that there exists
$(E,n,B)$ satisfying ($\epsilon\rightarrow 0$)
\begin{align}
\begin{split}\label{6.1}
&\sup\limits_{t\in[0,T]}(\|E^{\epsilon}-E\|_{H^{\tilde{s}+1}}+\|n^{\epsilon}-n\|_{H^{\tilde{s}}}
+\|B^{\epsilon}-B\|_{H^{\tilde{s}}\cap
\dot{H}^{-1}})\rightarrow 0,\\
&\sup\limits_{t\in[0,T]}(\|E^{\epsilon}_{t}-E_{t}\|_{H^{\tilde{s}-1}}+\|n^{\epsilon}_{t}-n_{t}\|_{H^{\tilde{s}-1}\cap\dot{H}^{-1}}
+\|B^{\epsilon}_{t}-B_{t}\|_{H^{\tilde{s}-2}\cap\dot{H}^{-2}})\rightarrow
0,
\end{split}
\end{align}
for all $\tilde{s}<s$. Moreover, for all $\tilde{s}<s$, we have
\begin{align*}
&n^{\epsilon}E^{\epsilon}\rightarrow nE\ \mathrm{in}\
C([0,T];H^{\tilde{s}}),\ E^{\epsilon}\times
B^{\epsilon}\rightarrow E\times B,\ \mathrm{in}\
C([0,T];H^{\tilde{s}}),\\
& |E|^{2} \rightarrow |E|^{2} \ \mathrm{in}\
C([0,T];H^{\tilde{s}+1}),\ E^{\epsilon}\times
\overline{E^{\epsilon}}\rightarrow E\times \overline{E},\
\mathrm{in}\ C([0,T];H^{\tilde{s}+1}).
\end{align*}
Now letting $\epsilon\rightarrow 0$ in \eqref{3.1a}-\eqref{3.1d},
and using the above strong convergence properties, we finally see
that $(E,n,B)$ is a solution of the original magnetic Zakharov
system \eqref{1.1}-\eqref{1.2}. Furthermore, by the boundedness
property \eqref{4.2}-\eqref{A4.1} and the strong convergence
property \eqref{6.1}, we have
\begin{equation*}
\begin{array}{cl}
(E,n,B)\in C\big([0,T];H^{s+1}\oplus H^{s}\oplus (H^{s }\cap
\dot{H}^{-1})\big),&\\
(E_{t},n_{t},B_{t})\in C\big([0,T];H^{s-1}\oplus (H^{s-1}\cap
\dot{H}^{-1})\oplus (H^{s-2}\cap \dot{H}^{-2})\big).&
\end{array}
\end{equation*}
From Proposition \ref{prop4.1}, the existence time $T$ depends on
the norm of the initial data. In fact, if $T_{\mathrm{max}}$ is
the maximal lifespan of the solution, then  either
$T_{\mathrm{max}}=\infty$ or $T_{\mathrm{max}}<\infty$ and
$$
\|E(t)\|_{H^{s+1}}+\|n(t)\|_{H^{s}}+\|n_{t}(t)\|_{H^{s-1}\cap
\dot{H}^{-1}}+\|B(t)\|_{H^{s}\cap
\dot{H}^{-1}}+\|B_{t}(t)\|_{H^{s-2}\cap
\dot{H}^{-2}}\rightarrow\infty
$$
as $t\rightarrow T_{\mathrm{max}}$. Hence, the local existence
part of Theorem \ref{thm1.1} is proved.

For the uniqueness of the system \eqref{1.1}-\eqref{1.2}, one can
apply the same methods used in Lemma \ref{lem5.1}, hence, the
proof of uniqueness is essentially the same as the proof of
\eqref{5.1}, and we omit the details. We thus finish the proof of
Theorem \ref{thm1.1}. \qed

\begin{remark}
Taking $\epsilon\rightarrow 0$ in \eqref{3.2} and \eqref{3.3}, and
using \eqref{4.2}-\eqref{A4.1}, \eqref{6.1}, then we see that the
solution of \eqref{1.1} also satisfies the conservation laws
\eqref{2.1}-\eqref{2.2}. Besides, the solution $(E,n,B)$ obtained
in Theorem \ref{1.1} depends continuously on the initial data in
the following sense: There exists $T_{1}>0$ depending on $R$ such
that if $E_{0}^{j}\rightarrow E_{0}$ in $H^{s+1}$,
$n_{0}^{j}\rightarrow n_{0}$ in $H^{s}$, $n_{1}^{j}\rightarrow
n_{1}$ in $H^{s-1}\cap \dot{H}^{-1}$, $B_{0}^{j}\rightarrow B_{0}$
in $H^{s}\cap \dot{H}^{-1}$, $B_{1}^{j}\rightarrow B_{1}$ in
$H^{s-2}\cap \dot{H}^{-2}$, and if $(E^{j},n^{j},B^{j})$ is the
corresponding solution of \eqref{1.1} with initial data
$(E_{0}^{j},n_{0}^{j},n_{1}^{j},B_{0}^{j},B_{1}^{j})$, then
$(E^{j},n^{j},B^{j})$ is defined on $[0,T_{1}]$ when $j$ is large,
and
\begin{align*}
&\sup\limits_{t\in[0,T_{1}]}(\|E^{j}-E\|_{H^{\tilde{s}+1}}+\|n^{j}-n\|_{H^{\tilde{s}}}
+\|B^{j}-B\|_{H^{\tilde{s}}\cap
\dot{H}^{-1}})\rightarrow 0,\\
&\sup\limits_{t\in[0,T_{1}]}(\|E^{j}_{t}-E_{t}\|_{H^{\tilde{s}-1}}+\|n^{j}_{t}-n_{t}\|_{H^{\tilde{s}-1}\cap\dot{H}^{-1}}
+\|B^{j}_{t}-B_{t}\|_{H^{\tilde{s}-2}\cap\dot{H}^{-2}})\rightarrow
0
\end{align*}
for all $\tilde{s}<s$ as $j\rightarrow \infty$.
\end{remark}

Note that Theorem \ref{thm1.1} needs the additional assumption
$n_{1}, B_{0} \in \dot{H}^{-1}$, $B_{1}\in \dot{H}^{-2}$. As
described in Section 1, this assumption is rather strong. In fact,
this condition can be removed by splitting the initial data into
low frequency part and high frequency part. Denote
$\varphi(\xi)\in C_{c}^{\infty}(\mathbb{R}^{d})$ such that $0\leq
\varphi(\xi)=\varphi(|\xi|)\leq 1$, $\varphi\equiv 1$ if
$|\xi|\leq 1$ and $\varphi\equiv 0$ if $|\xi|\geq 2$. For any
given $f\in H^{r}(\mathbb{R}^{d})$, $r\in \mathbb{R}$, we
decompose $f=f_{L}+f_{H}$, where
$$
\widehat{f_{L}}=\varphi(\xi)\hat{f},\
\widehat{f_{H}}=(1-\varphi(\xi))\hat{f}.
$$
So one can easily see that $f_{L}\in H^{k}$ for all $k\in
\mathbb{R}$ and $f_{H}\in \dot{H}^{l}\cap H^{l}$ for all $l\leq
r$. Furthermore, there holds
$$
\|f_{L}\|_{H^{k}}\leq \max\{5^{\frac{k-r}{2}},\}\|f\|_{H^{r}},\
\|f_{H}\|_{H^{l}\cap \dot{H}^{l}}\leq
\max\{2^{-\frac{r}{2}},1\}\|f\|_{H^{r}}.
$$
In this way, we can decompose $n_{1}\in H^{s}$ as
$n_{1}=n_{1L}+n_{1H}$ with $n_{1L}\in H^{k}$ for all $k\in
\mathbb{R}$, $n_{1H}\in H^{l}\cap\dot{H}^{l}$ for all $l\leq s$,
and in particular $n_{1H}\in H^{s}\cap\dot{H}^{-1}$. Moreover we
have
\begin{equation}\label{6.2}
\|n_{1L}\|_{H^{k}}\leq C(k,s)\|n_{1}\|_{H^{s}},\ \forall\ k\in
\mathbb{R},\ \|n_{1H}\|_{H^{s}\cap\dot{H}^{-1}}\leq
\|n_{1}\|_{H^{s}}.
\end{equation}
Similarly, for given $B_{0}\in H^{s}$, $B_{1}\in H^{s-2}$, we have
$B_{0}=B_{0L}+B_{0H}$, $B_{1}=B_{1L}+B_{1H}$, where $B_{0L},\
B_{1L}\in H^{k}$ ($k\in \mathbb{R}$), $B_{0H}\in
H^{s_{1}}\cap\dot{H}^{s_{1}}$ for all $s_{1}\leq s$, $\ B_{1H}\in
H^{s_{2}}\cap\dot{H}^{s_{2}}$ for all $s_{2}\leq s-2$ and
\begin{align}\label{6.3}
\begin{split}
&\|B_{0L}\|_{H^{k}}\leq C(k,s)\|B_{0}\|_{H^{s}},\
\|B_{1L}\|_{H^{k}}\leq C(k,s-2)\|B_{1}\|_{H^{s-2}},\ \forall\
k\in \mathbb{R},\\
 &\|B_{0H}\|_{H^{s}\cap\dot{H}^{-1}}\leq
\|B_{0}\|_{H^{s}},\ \|B_{1H}\|_{H^{s-2}\cap\dot{H}^{-2}}\leq
C(s-2) \|B_{0}\|_{H^{s-2}}.
\end{split}
\end{align}

Now we set
\begin{equation}\label{6.4}
\tilde{n}=n-tn_{1L},\ \tilde{B}=B-B_{0L}-tB_{1L},
\end{equation}
and consider the equation
\begin{equation}\label{6.5}
\left\{\!\!
\begin{array}{lc}
iE_{t}+\nabla(\nabla\cdot E)-\alpha\nabla\times(\nabla\times
E)-(\tilde{n}+tn_{1L})E+iE\times (\tilde{B}+B_{0L}+tB_{1L})=0,&\\
\tilde{n}_{tt}-\Delta \tilde{n}=\Delta |E|^{2}+t\Delta n_{1L},&\\
\tilde{B}_{tt}+\Delta ^{2}\tilde{B}-\Delta
\tilde{B}=-i\Delta^{2}(E\times
\bar{E})-\Delta^{2}(B_{0L}+tB_{1L})+\Delta
(B_{0L}+t B_{1L})&\\
\end{array}
\right.
\end{equation}
with initial data
\begin{equation}\label{6.6}
E(0,x)=E_{0},\
(\tilde{n}(0,x),\tilde{n}_{t}(0,x))=(n_{0},n_{1H}),\
(\tilde{B}(0,x),\tilde{B}_{t}(0,x))=(B_{0H},B_{1H}).
\end{equation}
Note that the initial data \eqref{6.6} satisfies the condition
$\tilde{n}_{t}(0)\in \dot{H}^{-1}$, $\tilde{B}(0)\in
\dot{H}^{-1}$, $\tilde{B}_{t}(0)\in \dot{H}^{-2}$. We also remark
that if $(E,n,B)$ solves \eqref{1.1}-\eqref{1.2}, then
$(E,\tilde{n},\tilde{B})$ defined by \eqref{6.4} solves
\eqref{6.5}-\eqref{6.6}, and vice versa.

For the regular solution of equation \eqref{6.4}, a similar
argument as in Proposition \ref{prop2.1} gives that
$\|E(t)\|_{L^{2}}=\|E_{0}\|_{L^{2}}$ and
\begin{align}
&\frac{d}{dt}\Big(\|\nabla\cdot
E(t)\|_{L^{2}}^{2}+\alpha\|\nabla\times
E(t)\|_{L^{2}}^{2}+\frac{1}{2}\|\tilde{n}(t)\|_{L^{2}}^{2}+\frac{1}{2}\|\Lambda^{-1}\tilde{n}_{t}(t)\|_{L^{2}}^{2}+
\frac{1}{2}\|\Lambda^{-2}\tilde{B}_{t}(t)\|_{L^{2}}^{2}\nonumber\\
&\quad+\frac{1}{2}\|\tilde{B}(t)\|_{L^{2}}^{2}+\frac{1}{2}\|\Lambda^{-1}\tilde{B}(t)\|_{L^{2}}^{2}
+\int_{\mathbb{R}^{d}}\tilde{n}(t)|E(t)|^{2}dx+\int_{\mathbb{R}^{d}}tn_{1L}|E(t)|^{2}dx\nonumber\\
&\quad+i\int_{\mathbb{R}^{d}}\big(E(t)\times
\overline{E(t)}\big)\cdot
\tilde{B}(t)dx+i\int_{\mathbb{R}^{d}}\big(E(t)\times
\overline{E(t)}\big)\cdot (B_{0L}+tB_{1L})dx\Big)\nonumber\\
&=\int_{\mathbb{R}^{d}}n_{1L}|E(t)|^{2}dx+i\int_{\mathbb{R}^{d}}B_{1L}\cdot(E(t)\times
\overline{E(t)})dx-t\int_{\mathbb{R}^{d}}\Lambda
n_{1L}\Lambda^{-1}\tilde{n}_{t}dx\nonumber\\
&\quad-\int_{\mathbb{R}^{d}}\Lambda^{2}(
B_{0L}+tB_{1L})\Lambda^{-2}\tilde{B}_{t}dx-\int_{\mathbb{R}^{d}}(
B_{0L}+tB_{1L})\Lambda^{-2}\tilde{B}_{t}dx.\label{6.7}
\end{align}
Using \eqref{6.2}, \eqref{6.3} and the fact
$\|E\|_{L^{2}}=\|E_{0}\|_{L^{2}}$, we have
$$
\mathrm{RHS\ of\ \eqref{6.7}}\leq
C+C(1+t)(\|\Lambda^{-1}\tilde{n}_{t}\|_{L^{2}}^{2}+\|\Lambda^{-2}\tilde{B}_{t}\|_{L^{2}}^{2})
$$
and
\begin{align*}
&\left|\int_{\mathbb{R}^{d}}tn_{1L}|E(t)|^{2}dx\right|\leq
t\|n_{1L}\|_{L^{\infty}}\|E_{0}\|_{L^{2}}^{2}\leq Ct,\\
&\left|i\int_{\mathbb{R}^{d}}\big(E(t)\times
\overline{E(t)}\big)\cdot (B_{0L}+tB_{1L})dx\right|\leq C(1+t).
\end{align*}
Integrating \eqref{6.7}, and applying the same method given in
Lemma \ref{lem2.2} and Gronwall's inequality, we can bound the
quantity
$$
\|E\|_{H^{1}}^{2}+\|\tilde{n}\|_{L^{2}}^{2}+\|\tilde{n}_{t}\|_{\dot{H}^{-1}}^{2}+\|\tilde{B}\|_{L^{2}\cap
\dot{H}^{-1}}^{2}+\|\tilde{B}_{t}\|_{\dot{H}^{-2}}^{2}
$$
by the norm of initial data \eqref{6.6}. If we return to our
original system, then we can obtain the following result.
\begin{lemma}\label{lem6.1}
 Assume $(E,n,B)$ is a sufficiently regular solution to the magnetic
Zakharov system \eqref{1.1}-\eqref{1.2}, and let
$2\|E_{0}\|_{L^{2}}^{2}<\|Q\|_{L^{2}}^{2}$ in the case $d=2$,
where $Q=Q(x)$ is the ground state solution of
\begin{equation*}
\Delta Q-Q+Q^{3}=0,\ x\in \mathbb{R}^{2},
\end{equation*}
and $\|E_{0}\|_{H^{1}}$ small when $d=3$. Then
\begin{equation*}
\|E\|_{H^{1}}^{2}+\|n\|_{L^{2}}^{2}+\|n_{t}\|_{H^{-1}}^{2}+\|B\|_{L^{2}}^{2}+\|B_{t}\|_{H^{-2}}^{2}\leq
C,
\end{equation*}
here the constant $C$ depends on $t$,
$\|E_{0}\|_{H^{1}},\|n_{0}\|_{L^{2}},
\|n_{1}\|_{H^{-1}},\|B_{0}\|_{L^{2}},\|B_{1}\|_{H^{-2}}$.
\end{lemma}

Again Lemma \ref{lem6.1} implies the existence of weak solution
for the magnetic system.

\begin{theorem}\label{thm6.1}
If $E_{0}\in H^{1}$, $(n_{0},n_{1})\in L^{2}\oplus H^{-1}$,
$(B_{0},B_{1})\in L^{2}\oplus H^{-2}$, and the initial data
satisfying $2\|E_{0}\|_{L^{2}}^{2}<\|Q\|_{L^{2}}^{2}$ in the case
$d=2$ and $\|E_{0}\|_{H^{1}}$ small in $d=3$, then there exists a
weak solution $(E,n,B)$ for the system \eqref{1.1} in the sense of
distributions such that
$$
E\in L^{\infty}_{loc}(\mathbb{R}^{+};H^{1}),\ (n,n_{t})\in
L^{\infty}_{loc}(\mathbb{R}^{+};L^{2}\oplus H^{-1}),\ (B,B_{t})\in
L^{\infty}_{loc}(\mathbb{R}^{+};L^{2}\oplus H^{-2}).
$$
\end{theorem}

Due to \eqref{6.2} and \eqref{6.3}, we see that the low frequency
part of $n_{1}$, $B_{0}$ and $B_{1}$ appearing in the equation
\eqref{6.5} can be well controlled. Therefore, one can follow the
same procedure as in Section 3-Section 5 and then get the
existence and uniqueness of solution for the equation \eqref{6.5}
with initial data \eqref{6.6}, which in turn leads to Theorem
\ref{1.2}. Since this process is much the same as the proof of
Theorem \ref{1.1}, the details are omitted. Hence, in this way,
Theorem \ref{thm1.2} is proved.

\begin{center}

\end{center}

\end{document}